\documentclass{amsart}

\usepackage{a4wide}
\usepackage{latexsym}
\usepackage{amsfonts}
\usepackage{amsmath}
\usepackage{euscript}
\usepackage[all]{xy}
\usepackage{amscd}
\usepackage{amsthm}
\usepackage{amssymb}
\usepackage{verbatim}

\usepackage[T1]{CJKutf8}

\newtheorem{theorem}{Theorem}[subsection]
\newtheorem{proposition}[theorem]{Proposition}
\newtheorem{lemma}[theorem]{Lemma}
\newtheorem{corollary}[theorem]{Corollary}

\theoremstyle{definition}
\newtheorem{definition}[theorem]{Definition}
\newtheorem{example}[theorem]{Example}

\theoremstyle{remark}
\newtheorem{remark}[theorem]{Remark}

\numberwithin{equation}{section}

\def\Ob{\mathop{\rm Ob}}
\def\lim{\mathop{\varprojlim}}

\def\Hom{{\mathop{\rm Hom}}}
\def\Ext{{\mathop{\rm Ext}}}
\def\k{\underline{k}}
\def\End{{\mathop{\rm End}}}
\def\Aut{{\mathop{\rm Aut}}}
\def\Res{\mathop{\rm Res}}
\def\res{{\mathop{\rm res}}}
\def\Mor{\mathop{\rm Mor}}

\def\Id{\mathop{\rm Id}}
\def\Is{\mathop{\rm Is}}
\def\A{\mathcal{A}}
\def\C{\mathcal{C}}
\def\c{\mathbb{C}}
\def\D{\mathcal{D}}
\def\K{{\mathop{\sf K}}}
\def\hocolim{{\rm hocolim}}
\def\E{\mathcal{E}}
\def\H{{\mathop{\rm H}}}
\def\s{\mathfrak{S}}
\def\S{\mathcal{S}}
\def\P{\mathcal{P}}
\def\B{\mathfrak{B}}
\def\m{\mathfrak{M}}
\def\n{\mathfrak{N}}
\def\p{\mathfrak{P}}
\def\hotimes{\hat{\otimes}}
\def\GP{\mathop{(G\propto\P)}}
\def\gp{\mathop{G\propto\P}}
\def\Spec{{\rm MaxSpec}}
\def\i{\mathfrak{i}}
\def\F{\mathcal{F}}
\def\q{\mathfrak{q}}
\def\e{\mathbb{E}}
\def\l{\mathfrak{L}}
\def\Q{\mathcal{Q}}
\def\HH{{\mathop{\rm HH}}}
\def\Rad{\mathop{\rm Rad}}

\def\pd{\mathop{\rm proj.dim}}
\def\id{\mathop{\rm inj.dim}}

\def\cm{\mathop{\underline{\rm CM}}}
\def\Proj{\mathop{\rm Proj}}
\def\L{\mathcal{L}}

\def\O{\mathcal{O}}

\title[Support varieties]
 {Support varieties for transporter category algebras} 

\author{Fei Xu}

\address{Departament de Matem\`atiques, Universitat Aut\`onoma de Barcelona, 08193 Bellaterra, Espanya.}

\email{xu@mat.uab.cat}


\thanks{The author \begin{CJK*}{UTF8}{}
\CJKtilde \CJKfamily{gbsn}徐 斐
\end{CJK*} was supported by
a Beatriu de Pin\'os research fellowship from the government of
Catalonia of Spain, as well as a Grant MTM2010-20692 ``Analisis local en grupos y espacios topologicos''
from the Ministry of Science and Innovation of Spain.}


\begin{document}
\maketitle

\begin{abstract}
Let $G$ be a finite group. Over any finite $G$-poset $\P$ we may
define a transporter category as the corresponding Grothendieck
construction. The classifying space of the transporter category is 
the Borel construction on the $G$-space $B\P$, while the $k$-category 
algebra of the transporter category is the (Gorenstein) skew group algebra on the 
$G$-algebra $k\P$. 

We introduce a support variety theory for the category
algebras of transporter categories. It extends Carlson's support
variety theory on group cohomology rings to equivariant cohomology rings. 
In the mean time it provides a class of (usually non selfinjective) algebras to which
Snashall-Solberg's (Hochschild) support variety theory applies.
Various properties will be developed. Particularly we establish a
Quillen stratification for modules.
\end{abstract}

\section{Introduction}

Let $G$ be a finite group and $\P$ a finite $G$-poset. Throughout
this paper, we assume $k$ is an algebraically closed field of
characteristic $p$, dividing the order of $G$. We are interested in
a finite category $G\propto\P$, which is the Grothendieck
construction on the $G$-poset $\P$ and which we will call a
\textit{transporter category} in this paper. When $G=\{e\}$ is
trivial, $\{e\}\propto\P\cong\P$ and when $\P=\bullet$ is trivial,
$G\propto\bullet\cong G$. A transporter category $\gp$ is the
algebraic or categorical predecessor of the Borel construction
$EG\times_G B\P$ on the finite $G$-CW-complex $B\P$, in the sense
that $B\GP\simeq EG\times_G B\P$. Our interests in transporter
categories come from the fact that the equivariant cohomology ring
$\H^*_G(B\P,k)=\H^*(EG\times_G B\P,k)$ is Noetherian. Through an
algebraic construction of the equivariant cohomology ring, we may
introduce in a natural way modules over this ring and hence extend
Carlson's support variety theory for finite group algebras to one
for finite transporter category algebras.

Let us recall some historical developments in support variety
theory. Suppose a $G$-space $X$ is compact or paracompact with
finite cohomological dimension. Quillen \cite{Q1, Q2} proved that
$\H^*_G(X)$ is Noetherian. Following his notation, we put $\H_G(X)$
to be $\H^*_G(X)$ if $p=2$ or $\H^{ev}_G(X)$, the even part of the
equivariant cohomology ring, if $p
> 3$. When $X=\bullet$ is just a point fixed by $G$, the equivariant
ring reduces to the group cohomology ring and we shall write
$\H^*_G=\H^*_G(\bullet)$ and $\H_G=\H_G(\bullet)$. Quillen's work
began with the observation that the graded ring $\H_G(X)$ is
commutative Noetherian. It enabled him to define an affine variety
$V_{G,X}$ as the maximal ideal spectrum $\Spec\H_G(X)$, and
subsequently described it in terms of $V_E=V_{E,\bullet}=\Spec\H_E$,
where $E$ runs over the set of all elementary abelian $p$-subgroups
of $G$. This is what we nowadays refer to as the Quillen
stratification. Restricting to the special case of $X=\bullet$,
based on the fact that $\Ext^*_{kG}(M,M)$ is finitely generated over
$\H^*_G\cong\Ext^*_{kG}(k,k)$, Carlson \cite{C} extended Quillen's
work by attaching to every finitely generated $kG$-module $M$ a
subvariety of $V_G=V_{G,\bullet}$, denoted by
$V_G(M)=\Spec\H_G/I_G(M)$, called the (cohomological) support
variety of $M$, where $I_G(M)$ is the kernel of the following map
$$
\phi_M = -\otimes_k M : \H^*_G\cong\Ext^*_{kG}(k,k) \to
\Ext^*_{kG}(M,M).
$$
Especially since $\phi_k$ is the identity, $V_G=V_G(k)$. Following
Carlson's construction, Avrunin and Scott \cite{AS} quickly
generalized the Quillen stratification from $V_G$ to $V_G(M)$. By
showing that support varieties are well-behaved with respect to
module operations, gradually Benson, Carlson and many others
developed a remarkable theory, being a significant progress in group
representations and cohomology. Since then, some other analogous
support variety theories have been introduced for restricted Lie
algebras \cite{FPa}, for finite group schemes \cite{BFS, FPe}, for
complete intersections \cite{AB} and for certain finite-dimensional
algebras \cite{L2, EHSST, SS}.

Quillen's work on equivariant cohomology rings has not been fully
exploited, partially because there existed no suitable modules which
$\H^*_G(X)$ (hence $\H_G(X)$) acts on or maps to, as in Carlson's
theory. In this article we attempt to use category algebras to solve
the problem: if $X=B\P$ comes from a finite $G$-poset, then we
consider the category algebra $k(G\propto\P)$ of the transporter
category $G\propto\P$, based on which we will generalize Carlson's
theory. In fact, let $\k$ be the trivial $k\GP$-module (see Section
2.2). Then $\Ext^*_{k\GP}(\k,\k)$ is a graded commutative ring and
there exists a natural ring isomorphism
$$
\Ext^*_{k\GP}(\k,\k)\cong\H^*(EG\times_G B\P,k)=\H^*_G(B\P,k).
$$
We shall call the above ring the \textit{ordinary cohomology ring}
of $k\GP$ (instead of the equivariant cohomology ring), as opposed
to the \textit{Hochschild cohomology ring} of $k\GP$. Then we define
$V_{\gp}=V_{G,B\P}=\Spec\H_G(B\P)$. The virtue of having an entirely
algebraic construction of the equivariant cohomology theory is that
it allows us to consider
$$
\Ext^*_{k\GP}(\m,\n)
$$
for any finitely generated $\m,\n\in k\GP$-mod, and moreover
construct a map
$$
\phi_{\m} = -\hotimes\m : \Ext^*_{k\GP}(\k,\k) \to
\Ext^*_{k\GP}(\m,\m).
$$
Here $\hotimes$ is the tensor product in the closed symmetric
monoidal category ($k\GP$-mod, $\hotimes, \k$). Note that $\k$
serves as the tensor identity. Since we have shown in \cite{X4} that
$\Ext^*_{k\GP}(\m,\n)$ is finitely generated over the ordinary
cohomology ring, we may define the support variety of $\m\in
k\GP$-mod as $V_{\gp}(\m)=\Spec\H_G(B\P)/I_{\gp}(\m)$, where
$I_{\gp}(\m)$ is the kernel of $\phi_{\m}$. Especially
$V_{\gp}=V_{\gp}(\k)$. When $\P=\bullet$, the is exactly Carlson's
construction because $k(G\propto\bullet)\cong kG$, $\k$ becomes the
trivial $kG$-module $k$ and $\hotimes$ reduces to $\otimes_k$ under
the circumstance.

As a surprising consequence of our investigations of transporter
category algebras, we assert that Snashall-Solberg's (Hochschild)
support variety theory (for Gorenstein algebras) applies to every block of a finite
transporter category algebra. Furthermore, our support variety
theory is closely related with Snashall-Solberg's as what happens in
the case of group algebras and their blocks. A notable point is that 
the block algebras of a transporter category algebra are usually non-selfinjective 
and non-commutative, opposing to the cases of (selfinjective) 
Hopf algebras \cite{C, BFS, FPe} and of commutative Gorenstein algebras \cite{AB} 
considered by others.

This paper is organized as follows. Section 2 recalls the
definitions of a transporter category, a category algebra and the
category cohomology. Various necessary constructions are recorded
for the convenience of the reader. Here we show a transporter
category algebra is Gorenstein and the ordinary cohomology ring of
such an algebra is identified with an equivariant cohomology ring.
Then in Section 3, we define the support variety for a module over a
transporter category algebra. To motivate the reader, we describe
the relevant works of Carlson, Linckelmann and Snashall-Solberg,
before we develop some standard properties. Sections 4 and 5 contain
various properties of support varieties, including the Quillen
stratification for modules, as well as results related with
sub-transporter categories and tensor products. Notably we propose a
concept of support varieties of Mackey functors.

\section{preliminaries}

In this section, we recall the definition of a transporter category
and some background in category algebras. Throughout this article we
will only consider finite categories, in the sense that they have
finitely many morphisms. Thus a group $G$, or a $G$-poset $\P$, is
always finite.

The morphisms in a poset are customarily denoted by $\le$. The
\textit{dimension} of a poset $\P$, $\dim\P$, is defined to be the
maximal integer $n$ such that a chain of non-isomorphisms $x_0 < x_1
< \cdots < x_n$ exists in $\P$.

Any $G$-set is usually regarded as a $G$-poset with trivial
relations. One the other hand, since in a $G$-poset $\P$, both
$\Ob\P$ and $\Mor\P$ are naturally $G$-sets, we shall use
terminologies for $G$-sets in our situation without further
comments.

\subsection{Transporter categories as Grothendieck constructions}

We deem a group as a category with one object, usually denoted by
$\bullet$. The identity of a group $G$ is always named $e$. We say a
poset $\P$ is a $G$-poset if there exists a functor $F$ from $G$ to
$\mathfrak{sCats}$, the category of small categories, such that
$F(\bullet)=\P$. It is equivalent to say that we have a group
homomorphism $G \to \Aut(\P)$. The \textit{Grothendieck
construction} on $F$ will be called a \textit{transporter category}.

\begin{definition} Let $G$ be a group and $\P$ a $G$-poset. The
transporter category $G\propto \P$ has the same objects as $\P$,
that is, $\Ob(G\propto \P)=\Ob\P$. For $x ,y \in \Ob(G\propto \P)$,
a morphism from $x$ to $y$ is a pair $(g, gx\le y)$ for some $g\in
G$.
\end{definition}

If $(g, gx\le y)$ and $(h,hy\le z)$ are two morphisms in
$G\propto\P$, then their composite is easily seen to be $(hg,
(hg)x\le z)=(h,hy\le z)\circ(g, gx\le y)$.

\begin{remark} One can check directly that if
$\Hom_{G\propto\P}(x,y)\ne\emptyset$ then both
$\Aut_{G\propto\P}(x)$ and $\Aut_{G\propto\P}(y)$ act freely on
$\Hom_{G\propto\P}(x,y)$. This simple observation is quite useful to
us.

We note that for each $x\in\Ob\GP=\Ob\P$, $\Aut_{\gp}(x)$ is exactly
the isotropy group of $x$. For the sake of simplicity, we will often
write $G_x=\Aut_{\gp}(x)$, and $[x]=Gx$, the orbit of $x$. Note that
$[x]$ is a $G$-set, consisting of exactly the objects in $\gp$ that
are isomorphic to $x$.

If $\P_1$ is a $G_1$-poset and $\P_2$ is a $G_2$-poset, then
$$
(G_1\times
G_2)\propto(\P_1\times\P_2)\cong(G_1\propto\P_1)\times(G_2\propto\P_2).
$$
We will utilize this in Section 5.2.
\end{remark}

From the definition one can easily see that there is a natural
embedding $\iota : \P\hookrightarrow G\propto\P$ via $(x \le y)
\mapsto (e, x\le y)$. On the other hand, the transporter category
admits a natural functor $\pi : G\propto\P \to G$, given by
$x\mapsto \bullet$ and $(g,gx\le y)\mapsto g$. Thus we always have a
sequence of functors
$$
\P {\buildrel{\iota}\over{\hookrightarrow}} G\propto\P
{\buildrel{\pi}\over{\longrightarrow}} G
$$
such that $\pi\circ\iota(\P)$ is the trivial subgroup or subcategory
of $G$. Topologically it is well known that $B(G\propto\P)\simeq
\hocolim_{G}B\P\simeq EG\times_G B\P$. Passing to classifying
spaces, we obtain a fibration sequence
$$
B\P {\buildrel{B\iota}\over{\longrightarrow}} EG\times_G B\P
{\buildrel{B\pi}\over{\longrightarrow}} BG .
$$
Forming the transporter category over a $G$-poset eliminates the
$G$-action, and thus is the algebraic analogy of introducing a Borel
construction over a $G$-space. This is the first instance explaining
why a transporter category has anything to do with the equivariant
cohomology theory.

\begin{example}
If $G$ acts trivially on $\P$, then $G\propto\P = G\times \P$. In
this case for any $x\in\Ob(G\times\P)$, $G_x=G$.
\end{example}

\begin{example} Let $G$ be a finite group and $H$ a subgroup. We consider
the set of left cosets $G/H$ which can be regarded as a $G$-poset:
$G$ acts via left multiplication. The transporter category $G\propto
(G/H)$ is a connected groupoid whose skeleton is isomorphic to $H$.
In this way one can recover all subgroups of $G$, up to category
equivalences. Making Grothendieck constructions on transitive
$G$-sets reveals the isotropy groups.

For an arbitrary $G$-poset $\P$ and $x\in\Ob\P$, we have a category
equivalence $G\propto[x]\simeq G_x$, see Remark 2.1.2. For instance,
$G$ acts on $\P=\S_p$, the poset of non-identity $p$-subgroup of
$G$, by conjugation. Then $G_x=N_G(x)$ for every $x\in
\Ob(G\propto\S_p)=\Ob\S_p$. The isotropy group of $x$, $G_x$, is
frequently identified with the transporter category $G_x\times x
\cong G_x\times\bullet$.
\end{example}

In the upcoming Section 2.2 we will see that a category equivalence
$\D\to\C$ induces a Morita equivalence between their category
algebras, $k\D\simeq k\C$, as well as a homotopy equivalence
$B\D\simeq B\C$ (see \cite{W}). It means there is no essential
difference between $H$ and $G\propto(G/H)$ as far as we concern.
Hence it makes sense if we deem transporter categories as
generalized subgroups for a fixed finite group.

\subsection{Category algebras and their representations}

We recall some facts about category algebras. The reader is referred
to \cite{W,X2} for further details. Let $\C$ be a finite category
and $k$ a field. One can define the \textit{category algebra} $k\C$,
which, as a vector space, has a basis the set of all morphisms in
$\C$. We then define a product on the base elements and extend it
linearly to $k\C$. The product $\alpha*\beta$ of any two base
elements $\alpha, \beta \in\Mor\C$ is defined to be $\alpha\beta$,
if they are composable, or zero otherwise. It is a
finite-dimensional associative algebra with identity
$1=\sum_{x\in\Ob\C}1_x$. The category algebra $k\C$ possesses a
distinguished module $\k$, called the \textit{trivial module}. It
can be defined as $\k=k\Ob\C$. If $\alpha$ is a base element of
$k\C$ and $x\in\Ob\C$, we ask $\alpha\cdot x= y$ if
$\alpha\in\Hom_{\C}(x,y)$. Otherwise we set $\alpha\cdot x=0$.

When $\C$ is a group, $k\C$ is exactly the group algebra and $\k=k$.
All modules we consider here are finitely generated left modules,
unless otherwise specified. As a convention, throughout this
article, the $kG$-modules are usually written as $M, N$ etc, while
the modules of a (non-group) category algebra $k\C$ are denoted by
$\m,\n$ etc., except some distinguished modules, namely $\k$ and, in
the special case of $\C=\gp$, $\kappa_M$ which are obtained from
$kG$-modules (to be defined shortly).

A $k$-\textit{representation} of $\C$ is a covariant functor from
$\C$ to $Vect_k$, the category of finite dimensional $k$-vector
spaces. All representations of $\C$ form the functor category
$Vect_k^{\C}$. By a theorem of B. Mitchell (see \cite{W}), the
finitely generated left $k\C$-modules are the same as the
$k$-representations of $\C$, in the sense that there exists a
natural equivalence
$$
Vect_k^{\C} \simeq k\C\mbox{-mod}.
$$
It is often helpful to utilize the underlying functor structure of a
module. For instance, upon the preceding category equivalence we can
alternatively define the trivial module $\k$ as a constant functor
taking $k$ as its value at every object of $\C$. Meanwhile since
$Vect_k$ is a symmetric monoidal category, $Vect_k^{\C}$ inherits
this structure. It means there exists an (internal) \textit{tensor
product}, or the pointwise tensor product, written as $\hotimes$,
such that for any two $k\C$-modules $\m,\n$,
$(\m\hotimes\n)(x):=\m(x)\otimes_k\n(x)$. Let $\alpha\in\Mor\C$ be a
base element of $k\C$. Then $\alpha$ acts on $\m\hotimes\n$ via
$\alpha\otimes\alpha$. Obviously $\k$ is the \textit{tensor
identity} with respect to $\hotimes$ and
$\m\hotimes\n\cong\n\hotimes\m$. The category $k\C$-mod has function
objects, also called \textit{internal hom} \cite{X2}, in the sense
that, for $\l, \m, \n \in k\C$-mod,
$$
\Hom_{k\C}(\l\hotimes\m,\n)\cong\Hom_{k\C}(\l,\mathcal{H}om(\m,\n)).
$$

We record the basic tools for comparing category algebras and their
modules. When $\tau : \D \to \C$ is a functor between two finite
categories, there are adjoint functors for comparing their
representations. The functor $\tau$ usually does not induce an
algebra homomorphism from $k\D$ to $k\C$. However it does induce an
exact functor, called the \textit{restriction} along $\tau$,
$$
{\Res}^{\C}_{\D} : k\C\mbox{-mod} \to k\D\mbox{-mod}.
$$
If we regard a $k\C$-module as a functor, then its restriction is
the precomposition with $\tau$. If we consider the functor $\pi :
G\propto\P \to G$, then any $kG$-module $M$ restricts to a
$k(G\propto\P)$-module, written as $\kappa_M=\Res^G_{\gp}M$, with
only one exception $\k=\kappa_k=\Res^G_{\gp}k$. It is known that
$\mathcal{H}om(\kappa_M,\kappa_N)\cong\kappa_{\Hom_k(M,N)}$ for any
two $M, N\in kG$-mod. In this paper, for the sake of simplicity, if
$\D\to\C$ is a functor and $\m$ is a $k\C$-module, when it will not
cause confusions, we sometimes neglect $\Res^{\C}_{\D}$ and regard
$\m$ (instead of $\Res^{\C}_{\D}\m$) as a $k\D$-module.

The functor $\Res^{\C}_{\D}$ is equipped with two adjoints: the
\textit{left and right Kan extensions} along $\tau$
$$
LK^{\C}_{\D}, RK^{\C}_{\D} : k\D\mbox{-mod} \to k\C\mbox{-mod}.
$$
The definition of the left and right Kan extensions depend on the
so-called over-categories and under-categories, respectively.
Despite their seemingly abstract definitions, they are quite
computable and thus play an important role in category
representations and cohomology, see \cite{X1, X3, X4}, as well as
Sections 2.5, 3.4 and 4.1.

Note that our notations for the restriction and Kan extensions are
slightly different from earlier articles. The reason is that in this
place we feel it is necessary to emphasize the two categories
involved in order to make the notations more indicative.

\subsection{EI categories} When a category $\C$ is an \textit{EI category}, that is, every
endomorphism is an isomorphism, there exists a natural partial order
on the set of isomorphism classes of objects in $\C$ \cite{W}.
Indeed the partial order is given by $[x]\le[y]$ if and only if
$\Hom_{\C}(x,y)\ne\emptyset$. Groups, posets and transporter
categories are all EI categories. The upshot is that this partial
order allows us to give a filtration of any $k\C$-module $\m$. Let
us call $x\in\Ob\C$ an $\m$-\textit{object} if $\m(x)\ne 0$. Assume
$x$ is minimal among all $\m$-object, then we can define a submodule
$\m_{\hat x}$ such that $\m_{\hat x}(y)=\m(y)$ unless $y\cong x$ in
$\C$, in which case $\m_{\hat x}(y)=0$. Then we have a short exact
sequence
$$
0 \to \m_{\hat x} \to \m \to \m/\m_{\hat x} \to 0.
$$
Denote by $\m_x=\m/\m_{\hat x}$. This is an example of the so-called
atomic modules of $k\C$. A $k\C$-module $\m$ is called
\textit{atomic} if as a functor $\m$ takes zero values on all but
one isomorphism classes of objects. From the above analysis we see
that every $\m$ admits a filtration with atomic modules as
composition factors. Obviously from any module $\m$ and an
$\m$-object $x$ we may define an atomic module by brutal truncation
(the restriction along the inclusion $[x]\to\C$). Abusing
terminology, we always write such modules as $\m_x$.

Assume $\C$ is finite EI. Then we can characterize projective and
injective $k\C$-modules. Recall from \cite{W} that each
indecomposable left (resp. right) projective $k\C$-module, up to
isomorphism, is a direct summand of $k\Hom_{\C}(x,-)=k\C\cdot 1_x$
(resp. $k\Hom_{\C}(-,x)=1_x\cdot k\C$) for some $x\in\Ob\C$. The
(indecomposable) right (resp. left) injective modules can be
obtained as $k$-duality.

\begin{lemma}
\begin{enumerate}
\item If $\p$ is a projective (or injective) $k\GP$-module, then $\p_x$ is a
projective $k(G\propto[x])$-module, for all $\p$-objects $x$.

\item A $k\GP$-module $\m$ is of finite projective (equivalently, finite injective)
dimension if and only if $\m_x$ is a projective
$k(G\propto[x])$-module, for all $\m$-objects $x$. Under the
circumstance, both $\pd\m$ and $\id\m$ are bounded by $\dim\P$.

\item The transporter category algebra $k\GP$ is Gorenstein, which
means as either the left or the right regular module, it has finite injective
dimension.
\end{enumerate}

\begin{proof} These are direct consequences of the fact that if
$\Hom_{\gp}(x,y)\ne\emptyset$, then $k\Hom_{\gp}(x,y)$ is both free
$kG_x$- and $kG_y$-module, along with the characterizations of
projectives and injectives. We only prove (2).

Let $\p_* \to \m \to 0$ be a finite projective resolution. Then its
restriction to $G\propto[x]$ is a finite projective resolution of
$\m_x$. Since $G\propto[x]\simeq G_x$, their category algebras are
Morita equivalent and thus $k(G\propto[x])$ is selfinjective. Then
the finite projective resolution $(\p_x)_*\to\m_x\to 0$ splits and
$\m_x$ is a projective (and injective) $k(G\propto[x])$-module. On
the other hand, assume $\m$ satisfies the property that $\m_x$ is
projective (or zero) for every object $x$. Take the projective cover
$\p \to \m$. We immediately know that $\p_x\to\m_x$ splits as a
$k(G\propto[x])$-map. Let $[y_1],\cdots,[y_n]$ be the isomorphism
classes of all minimal objects among $\m$-objects. Then
$\bigoplus_{i=1}n\p_{y_i}$ is the projective cover of $\m_y$ which
implies $\bigoplus_{i=1}^n\p_{y_i}\cong\m_y$. Hence if we examine
the kernel $\m'$ of $\p\to\m$, it has the property that $\m'_y$ is a
projective $k(G\propto[y])$-module for all objects $y$ such that
$[y]>[y_i]$ for some $1\le i \le n$, and zero else. Repeat the same
process for $\m'$, eventually we will obtain the finite projective
resolution of $\m$. It implies $\pd\m\le\dim\P$.

As to the injective dimension, we consider the right $k\GP$-module
$\m^{\wedge}$. It satisfies $\m^{\wedge}(x)=\m(x)^{\wedge}$. Similar
to the left module situation, it has a finite projective resolution.
When we dualize it, it becomes an injective resolution of $\m$.
\end{proof}
\end{lemma}

Since $G\propto[x]\simeq G_x$ (see Example 2.1.4.), in the first two
statements, we may replace $\p_x$, $\m_x$, $k(G\propto[x])$ by
$\p(x)$, $\m(x)$, $kG_x$, respectively.

It is helpful to give the following characterization of a transporter category algebra 
as a skew group algebra. Recall from \cite[Chapt. III, Section 4]{ARS} that if a $k$-algebra 
$A$ is a $G$-algebra, then we may define the \textit{skew 
group algebra} $A[G]$ to be the $k$-vector space $A\otimes_k kG$ equipped with 
a multiplication rule determined by
$$
(a_1\otimes g_1)*(a_2\otimes g_2)=a_1(g_1\cdot a_2)\otimes g_1g_2,
$$
where $a_1, a_2 \in A$ and $g_1, g_2 \in G$. Here $g_1\cdot a_2$ means the image of $a_2$ under 
the action of $g_1$. For instance we immediately have $kG\cong (k\bullet)[G]$. The reader is referred 
to \cite{AR} and \cite{ARS} for further details 
and known results about skew group algebras. 
Proposition 2.2 of \cite{AR} says that $A[G]$ is Gorenstein if and only if $A$ is. Since it is easy to verify that 
$k\P$ is Gorenstein, with the following result we have another proof of $k\GP$ being Gorenstein.

\begin{lemma} There is an algebra isomorphism $k\GP\cong k\P[G]$.

\begin{proof} The isomorphism is given by $(gx\le y, g)\mapsto (gx\le y)\otimes g$ with inverse 
$(x\le y)\otimes h\mapsto (h(h^{-1}x)\le y,h)$ for $x, y \in\Ob\P$ and $g,h \in G$.
\end{proof}
\end{lemma}

The modular representation theory of $k\GP$ will be studied in another place.

\subsection{Category cohomology and spectrum} For any two $k\C$-modules it
makes sense to consider $\Ext^*_{k\C}(\m,\n)=\bigoplus_{i \ge
0}\Ext^i_{k\C}(\m,\n)$. Furthermore if $\m'$ and $\n'$ are two other
modules, the tensor product $\hotimes$ induces a \textit{cup
product} as follows
$$
\cup : \Ext^i_{k\C}(\m,\n)\otimes\Ext^j_{k\C}(\m',\n') \to
\Ext^{i+j}_{k\C}(\m\hotimes\m',\n\hotimes\n').
$$
In particular $\Ext^*_{k\C}(\k,\k)$ is a graded commutative ring and
we have a natural isomorphism $\Ext^*_{k\C}(\k,\k)\cong\H^*(B\C,k)$
\cite{X2}. This ring is called the \textit{ordinary cohomology ring}
of $k\C$ and it acts on $\Ext^*_{k\C}(\m,\n)$ via the cup product.
For any $k\C$-module $\m$, the Yoneda splice provides a ring
structure on $\Ext^*_{k\C}(\m,\m)$. When $\m=\k$, the cup product
and Yoneda splice give the same ring structure on
$\Ext^*_{k\C}(\k,\k)$. There exists a ring homomorphism whose image
lies in the center of the graded ring $\Ext^*_{k\C}(\m,\m)$
$$
-\hotimes\m : \Ext^*_{k\C}(\k,\k) \to \Ext^*_{k\C}(\m,\m).
$$
Moreover given a short exact sequence $0\to\m_1\to\m_2\to\m_3\to 0$
the resulting connecting homomorphism is a morphism of
$\Ext^*_{k\C}(\k,\k)$-modules.

In summary, finite category cohomology behaves very much like the
special case of finite group cohomology, except the finite
generation property. The ordinary cohomology ring of a category
algebra is usually far from finitely generated, but it is so when
$\C=G\propto\P$ is finite as it is isomorphic to the equivariant
cohomology ring $\H^*_G(B\P,k)$, see Section 2.5 and \cite{Q1, X3}.

The functor $\Res^{\C}_{\D}$ introduced earlier leads to a
restriction on cohomology
$$
\res^{\C}_{\D} : \Ext^*_{k\C}(\k,\k) \to \Ext^*_{k\D}(\k,\k).
$$
It coincides with the restriction $\H^*(B\C,k)\to\H^*(B\D,k)$,
induced by the continuous map $B\tau : B\D \to B\C$ between two
classifying spaces, see \cite{X2}.

From now on, we assume $k$ is algebraically closed. Throughout this
paper let us denote by
$$
\H(\C) = \Ext_{k\C}(\k,\k) = \left\{
\begin{array}{ll}
              \Ext^*_{k\C}(\k,\k) ,& \mbox{if the characteristic of}\ k= 2;\\
              \Ext^{2*}_{k\C}(\k,\k)\   ,& \mbox{if the characteristic
              of}\ k> 2.
             \end{array}
      \right.
$$
This graded ring is commutative. If $\Ext^*_{k\C}(\k,\k)$ is
Noetherian, then we can consider the {\it maximal ideal spectrum},
an affine variety,
$$
V_{\C}=\Spec\H(\C).
$$
Under the circumstance we will call $V_{\C}$ the \textit{variety} of
$\C$.

Assume both $\Ext^*_{k\C}(\k,\k)$ and $\Ext^*_{k\D}(\k,\k)$ are
Noetherian and there exists a functor $\tau : \D \to \C$. Since the
preimage of a maximal ideal is still a maximal ideal, there exists a
map between two varieties
$$
\iota^{\C}_{\D}:=(\res^{\C}_{\D})^{-1} : V_{\D} \to V_{\C}.
$$

These varieties and their subvarieties are our main subjects and
thus it is helpful if we can handle the restriction. In various
interesting cases the map is well understood. As an application of
Lemma 2.3.1 we obtain the following result.

\begin{lemma} Consider a transporter category $\gp$ and an object
$x\in\Ob\GP$. Then the inclusions $G_x\times x \to G\propto[x] \to
\gp$ induce two restrictions which fit into the following
commutative diagram
$$
\xymatrix{\Ext^*_{k\GP}(\m,\n) \ar[rr]^{\res^{\gp}_{G\propto[x]}}
\ar@{=}[d] &&
\Ext^*_{kG\propto[x]}(\m_x,\n_x) \ar[d]^{\cong}\\
\Ext^*_{k\GP}(\m,\n) \ar[rr]_{\res^{\gp}_{G_x\times x}} &&
\Ext^*_{kG_x}(\m(x),\n(x))},
$$
for any two modules $\m,\n\in k\GP$-mod.

\begin{proof} Let $\p \to \m \to 0$ be a projective resolution. Then
$\p_x \to \m_x \to 0$ remains a projective resolution of the
$kG\propto[x]$-module $\m_x$, by Lemma 2.3.1 (1). Hence the functor
$$
\res^{\gp}_{G\propto[x]} : k\GP\mbox{-mod} \to
k(G\propto[x])\mbox{-mod}
$$
induces a restriction
$$
\res^{\gp}_{G\propto[x]} : \Ext^*_{k\GP}(\m,\n) \to
\Ext^*_{k(G\propto[x])}(\m_x,\n_x).
$$
Similarly we have a map $\res^{\gp}_{G_x\times x}$. The commutative
diagram follows directly from the natural isomorphism
$\Ext^*_{k(G\propto[x])}(\m_x,\n_x)\cong\Ext^*_{k
G_x}(\m(x),\n(x))$.
\end{proof}
\end{lemma}

Because the two restrictions $\res^{\gp}_{G_x\times x} :
\Ext^*_{k\GP}(\m,\m)\to\Ext^*_{kG_x}(\m(x),\m(x))$ and
$\res^{\gp}_{G\propto[x]} :
\Ext^*_{k\GP}(\m,\m)\to\Ext^*_{k(G\propto[x])}(\m_x,\m_x)$ are ring
homomorphisms with respect to the Yoneda splice, there are two maps
$$
\iota^{\gp}_{G_x\times x} : V_{G_x\times x} \to V_{\gp},
$$
induced by $\res_x$ for $\m=\n=\k$, and similarly
$$
\iota^{\gp}_{G\propto[x]} : V_{G\propto[x]} \to V_{\gp}.
$$
Note that $V_{G_x\times x}=V_{G\propto[x]}$.

\subsection{Bar resolution and equivariant cohomology} One concept that
we will refer to is the \textit{bar resolution} $\B^{\C}_*$, a
combinatorially constructed projective resolution, of $\k \in
k\C$-mod. Given a functor $\tau : \D\to\C$, one can define
$\mathbb{C}_*(\tau/-,k)$, a complex of projective $k\C$-modules,
such that, for each $x\in\Ob\C$, $\tau/x$ is a finite category (the
category over $x$ or just an \textit{overcategory}) and
$\mathbb{C}_*(\tau/x,k)$ is the (normalized) chain complex resulting
from the simplicial $k$-vector spaces coming from the nerve of
$\tau/x$ \cite{X1}. Since overcategories are used in the proof of a
couple of results in the next section, we recall its definition. The
objects of $\tau/x$ are pairs $(d,\alpha)$ where $d\in\Ob\D$ and
$\alpha : \tau(d)\to x$ is a morphism in $\C$. A morphism from
$(d,\alpha)$ to $(d',\alpha')$ is a morphism $f : d \to d'$ in $\D$,
satisfying $\alpha=\tau(f)\alpha'$. For each finite category $\D$ we
define $\B^{\D}_*=\mathbb{C}_*(\Id_{\D}/-,k)$. It is known that
$$
LK^{\C}_{\D}\B^{\D}_*\cong\mathbb{C}_*(\tau/-,k).
$$
As an example, for a group $G$ the only overcategory $\Id_G/\bullet$
is the Cayley graph and thus $\B^G_*$ is the bar resolution of $k$
in group cohomology. For the convenience of the reader we recall
from \cite{X3} how we obtain equivariant cohomology from category
cohomology. Since we can explicitly calculate the unique
overcategory as $\pi/\bullet\cong\Id_G/\bullet\times\P$, it follows
that
$$
LK^G_{\gp}\B^{\gp}_*\cong\mathbb{C}_*(\pi/\bullet,k)\simeq\B^G_*\otimes\c_*(\P,k)
$$
Thus we have the following chain isomorphisms and homotopy, for any
$M\in kG$-mod,
$$
\begin{array}{ll}
\Hom_{k\GP}(\B^{\gp}_*,\kappa_M)&\cong\Hom_{kG}(LK^G_{\gp}\B^{\gp}_*,M)\\
&\simeq\Hom_{kG}(\B^G_*\otimes_k\mathbb{C}_*(\P,k),M)\\
&\cong\Hom_{kG}(\B^G_*,\Hom_k(\mathbb{C}_*(\P,k),M)).
\end{array}
$$
Note that $\c_*(\P,k)$ is a finite complex and consists in each
dimension of a permutation $kG$-module. From \cite[VII.7]{Br} we see
immediately $\Ext^*_{k\GP}(\k,\kappa_M)\cong\H^*_G(B\P,M)$ for any
$M\in kG$-mod.

In \cite{X3} we also established the following isomorphism
$$
\Ext^*_{k\GP}(\kappa_M,\n)\cong\Ext^*_{k\GP}(\k,
\mathcal{H}om(\kappa_M,\n)).
$$
Since $\mathcal{H}om(\kappa_M,\kappa_N)\cong\kappa_{\Hom_k(M,N)}$,
the $\Ext^*_{k\GP}(\k,\k)$-action on
$$
\Ext^*_{k\GP}(\kappa_M,\kappa_M)\cong\Ext^*_{k\GP}(\k,\mathcal{H}om(\kappa_M,\kappa_M))
$$
can be obtained from the canonical map $\k \to
\mathcal{H}om(\kappa_M,\kappa_M)\cong\kappa_{\End_k(M)}$, induced by
$1_k \mapsto \Id_M$. This is analogous to the group case.

\section{Support varieties for modules}

Let $k$ be an algebraically closed field. Suppose $\C=\gp$ is a
transporter category. We have learned that $\Ext^*_{k\GP}(\k,\k)$ is
Noetherian, over which $\Ext^*_{k\GP}(\m,\n)$ is a finitely
generated module. The action factors through the natural ring
homomorphism
$$
-\hotimes\m : \Ext^*_{k\GP}(\k,\k) \to \Ext^*_{k\GP}(\m,\m),
$$
and
$$
-\hotimes\n : \Ext^*_{k\GP}(\k,\k) \to \Ext^*_{k\GP}(\n,\n).
$$
Based on these, we are about to develop a support variety theory.
Since a transporter category $\gp$ is intimately related with $G$,
we will see our theory is a generalization of Carlson's theory.
Standard references for Carlson's theory include \cite[Chapter
5]{B2} and \cite[Chapters 8, 9, 10]{Evens}.

\subsection{Basic definitions}

For convenience, we shall assume $\gp$ is \textit{connected}, unless
otherwise specified. It is equivalent to saying that $\k\in
k\GP$-mod is indecomposable or that
$$
\Ext^0_{k\GP}(\k,\k)=\H^0(B\GP,k)=k.
$$
It does not mean that $\P$ is connected, see Example 2.1.4.

\begin{definition}
Given a transporter category $\gp$ and modules $\m,\n\in k\GP$-mod,
we write $I_{\gp}(\m)$ for the kernel of the map
$$
-\hotimes\m : \Ext_{k\GP}(\k,\k) \to \Ext_{k\GP}(\m,\m),
$$
and $V_{\gp}(\m)$, the \textit{support variety} of $\m$, for the
subvariety $\Spec(\H(\gp)/I_{\gp}(\m))$ of $V_{\gp}=V_{\gp}(\k)$.

Since both $\Ext^*_{k\GP}(\m,\m)$ and $\Ext^*_{k\GP}(\n,\n)$ act on
$\Ext^*_{k\GP}(\m,\n)$ via Yoneda splice, we further define
$I_{\gp}(\m,\n)$ as the annihilator of $\Ext_{k\GP}(\k,\k)$ on
$\Ext^*_{k\GP}(\m,\n)$. Then we set
$V_{\gp}(\m,\n)=\Spec(\H(\gp)/I_{\gp}(\m,\n))$.

We say a subvariety of $V_{\gp}$ is {\it trivial} if it is
$\mathfrak{m}=\Ext^+_{k\GP}(\k,\k)$, the positive part of
$\Ext_{k\GP}(\k,\k)$.
\end{definition}

Since $I_{\gp}(\m,\m)=I_{\gp}(\m)$, we get that
$V_{\gp}(\m,\m)=V_{\gp}(\m)$, and that
$$
I_{\gp}(\m,\n)\subset I_{\gp}(\m)+ I_{\gp}(\n).
$$
The latter implies
$$
V_{\gp}(\m,\n)\subset V_{\gp}(\m)\cap V_{\gp}(\n).
$$

Let $\P$ be a $G$-poset and $\Q$ an $H$-poset. Suppose there exists
a group homomorphism $\phi : H \to G$ as well as a functor $\theta :
\Q \to \P$ such that $\phi(h)\theta=\theta\circ h$ for all $h\in H$.
For convenience we record such a map as $(\phi, \theta) : (H,\Q)\to
(G,\P)$. They induce a functor
$$
\Theta : H\propto\Q \to \gp,
$$
which in turn gives rise to a restriction map
$$
\res^{\gp}_{H\propto\Q} : \Ext^*_{k\GP}(\k,\k) \to
\Ext^*_{k(H\propto\Q)}(\k,\k)
$$
and a map between varieties
$$
\iota^{\gp}_{H\propto\Q} : V_{H\propto\Q} \to V_{\gp}.
$$
For instance, Lemma 2.4.1 dealt with $(i,i) : (G_x,x) \to (G,\P)$
and $(\Id_G,i) : (G, [x]) \to (G,\P)$, where $i$ stands for the
inclusion.

\begin{example} Suppose $G$ acts trivially on $\P$. Then $(\Id_G, pt) : (G,\P)\to (G,\bullet)$
induce $\pi : \gp \to G$ and
$$
\res^{G}_{G\times\P} : \Ext^*_{kG}(k,k) \to
\Ext^*_{k(G\times\P)}(\k,\k)
$$
Since
$\Ext^*_{k(G\times\P)}(\k,\k)\cong\Ext^*_{kG}(k,k)\otimes\Ext^*_{k\P}(\k,\k)$
by the K\"unneth formula, and $\Ext^*_{k\P}(\k,\k)=\H^*(B\P,k)$ is
finite-dimensional, $V_{G\times\P}=\prod_n V_G$, where $n$ is the
number of connected components of $\P$. The restriction induces a
natural map $V_{G\times\P}\to V_G$. With the assumption that $G\propto\P=G\times\P$ is connected,
we actually have $n=1$ because $G\times\P$ is connected if and only if $\P$ is.
\end{example}

At this point, it seems to be a good idea to compare our theory with
those of Carlson, Linckelmann and Snashall-Solberg. By putting our
approach into the right context, we can better understand the ideas
and see what properties we may expect. Afterwards, we will present
various results concerning support varieties.

\subsection{Carlson's theory} When $\P=\bullet$ our theory
is just the theory of Carlson. However, combining recent works in
group and category cohomology, Carlson's theory can be recovered in
a more subtle way. To be explicit, if $\kappa_M\in k\GP$-mod for
some $M\in kG$-mod, then we have a commutative diagram \cite{X4}
$$
\xymatrix{\Ext^*_{kG}(k,k) \ar[rr]^{-\otimes M} \ar[d]_{\res^G_{G\propto\P}} &&
\Ext^*_{kG}(M,M) \ar[d]^{\res^G_{G\propto\P}}\\
\Ext^*_{k\GP}(\k,\k) \ar[rr]_(.5){-\hotimes\kappa_M} &&
\Ext^*_{k\GP}(\kappa_M,\kappa_M)}
$$
One can quickly deduce that the restriction map induces
$$
\Ext^*_{kG}(k,k)/I_G(M) \to \Ext^*_{k\GP}(\k,\k)/I_{\gp}(\kappa_M),
$$
and hence a finite (usually not surjective) map
$$
V_{\gp}(\kappa_M) \to V_G(M).
$$

If the Euler characteristic $\chi(\P,k)$ is invertible in $k$, then
by using the Becker-Gottlieb transfer map \cite{X4}, both vertical
maps are injective. Furthermore, if we let $\P=\S_p$, the poset of
non-identity $p$-subgroups of $G$, the left $\res^G_{G\propto\S_p}$
becomes an algebra isomorphism (see \cite[Chap. X, Section 7]{Br} and \cite[Part I, Sections 7 and
8]{DH} for instance). Hence we get
$$
I_G(M)\cong I_{G\propto\S_p}(\kappa_M)\ \ \mbox{and}\ \ V_G(M)\cong
V_{G\propto\S_p}(\kappa_M).
$$
It means that various properties of $V_G(M)$ can be rewritten for
$V_{G\propto\S_p}(\kappa_M)$. As an example we have a tensor product
formula
$$
V_{G\propto\S_p}(\kappa_M\hotimes\kappa_N)=V_{G\propto\S_p}(\kappa_{M\otimes
N})=V_G(M\otimes N)=V_G(M)\cap V_G(N)=V_{G\propto\S_p}(\kappa_M)\cap
V_{G\propto\P}(\kappa_N).
$$
Here the third equality comes from \cite[Theorem 5.7.1]{B2}. One can
similarly deduce other properties for $V_{G\propto\S_p}(\kappa_M)$
but we shall leave it to the reader as they are just reformulations
of known results for $V_G(M)$. Our interests really lie in
$V_{G\propto\S_p}(\m)$, or more generally $V_{\gp}(\m)$, for modules
$\m\ne\kappa_M$ for any $M\in kG$-mod.

In the above arguments, $\S_p$ may be replaced by various
$G$-subposets which are $G$-homotopy equivalent to it, see
\cite[Section 6.6]{B2}. A typical example is $\E_p$, the $G$-poset
of all elementary abelian $p$-subgroups of $G$. However we want to
emphasize that most of our results are established without
specifying a poset, see for instance Sections 4 and 5.

\subsection{Varieties in blocks}

A group algebra can be written as a direct product of
(indecomposable) block algebras $\prod_i b_i$. (Here for convenience
we denote by $b$, instead of $kGb$, a block algebra.) Each
indecomposable $kG$-module belongs to exactly one of these blocks in
the sense that all but one block act as zero on it. The block that
$k$ belongs to is called the principal block, denoted by $b_0$.

In \cite{L1}, Linckelmann introduced to each block algebra $b$ a
Noetherian graded commutative ring $\H^*(b)$, called the
\textit{block cohomology ring}. Then he showed that there exists a
natural injective homomorphism
$$
\H^*(b) \to \HH^*(b),
$$
and thus $\H^*(b)$ acts on $\Ext^*_b(M,M)$ via
$\HH^*(b)=\Ext^*_{b^e}(b,b)$ (this action is explained in Section
3.4), if $M\in b$-mod. Particularly if $b_0$ is the principal block
of a group algebra $kG$, then $\H^*(b_0)$ is isomorphic to
$H^*(G,k)$ and the above injection coincides with the composite of
two canonical maps
$$
\H^*(b_0)\cong\H^*(G,k)\to\HH^*(kG)\to\HH^*(b_0).
$$
Based on these, he was able to define support varieties for modules
of a block algebra as above in a natural way \cite{L2}, being a
refinement of Carlson's theory. Most significantly Linckelmann's
work brought Hochschild cohomology into the theory of support
varieties, which was taken up by Snashall and Solberg to develop a
new support variety theory using Hochschild cohomology rings, see
Section 3.4. Recently \cite{L3} Linckelmann had demonstrated that,
for a block algebra $b$, the inclusion $\H^*(b)\to\HH^*(b)$ induces
an isomorphism upon quotient out nilpotent elements. It implies that
the two support variety theories, of Lincklemann based on the block
cohomology ring and of Snashall-Solberg defined over the Hochschild
cohomology ring of a block of a finite group, are identical. See
\cite{BL} as well.

\subsection{Snashall-Solberg's theory}
Snashall and Solberg \cite{SS} developed a support variety theory
for certain finite-dimensional algebras using Hochschild cohomology
rings. Let $A$ be a finite-dimensional algebra, and $M, N$ two
finitely generated $A$-modules. Then there exists a natural action
of the Hochschild cohomology ring on Ext groups via the following
homomorphism
$$
\phi_M=-\otimes_A M : \Ext^*_{A^e}(A,A) \to \Ext^*_A(M,M).
$$
Based on Yoneda splice, one can introduce an action on
$\Ext^*_A(M,N)$ for any two $A$-modules.

For technical reasons, now let $A$ be an \textit{indecomposable}
algebra. Consequently $Z(A)$ becomes a commutative local algebra.
Let $\HH(A)=\Ext_{A^e}(A,A)$ be defined in the same fashion as
$\H(\C)$ in Section 2.4. Assume

\begin{enumerate}
\item[(Fg.1)] there is a graded Noetherian subalgebra $\H \subset \Ext_{A^e}(A,A)$
with $\H^0=\Ext^0_{A^e}(A,A)=Z(A)$; and

\item[(Fg.2)] for any $M,N\in A$-mod, $\Ext^*_A(M,N)$ is finitely generated over $\H$.
\end{enumerate}

Under the above assumptions, Snashall-Solberg considered the maximal
ideal spectrum $V_{\H}=\Spec\H$. Since $\H$ acts on $\Ext^*_A(M, N)$
for any two $A$-modules, assuming $I_{\H}(M,N)$ is the annihilator
they defined a subvariety by
$$
V_{\H}(M,N)=\Spec(\H/I_{\H}(M,N)).
$$
Write $I_{\H}(M)=I_{\H}(M,M)$. Then the (Hochschild) support variety
of $M \in A$-mod is given by
$$
V_{\H}(M) = \Spec(\H/I_{\H}(M)).
$$
They showed that $V_{\H}(M)=V_{\H}(M,A/\Rad A)=V_{\H}(A/\Rad A,M)$
and $V_{\H}(A/\Rad A)=V_{\H}$. Various satisfactory properties were
obtained in \cite{SS, EHSST}, see Theorem 3.4.4 below for a summary.
For future reference, we record the following definition.

\begin{definition}
A subvariety of $V_{\H}$ is called \textit{trivial} if it is
$\langle \Rad\H^0,\H^+\rangle$, where $\H^+$ consists of all the
positive degree elements in $\H$.
\end{definition}

Unfortunately the above conditions (Fg.1) and (Fg.2) are not met by
all finite-dimensional algebras, see \cite{X1}. Indeed they put
restrictions on the algebras that one may consider. For example, two
necessary conditions are that the algebra $A$ has to be Gorenstein
\cite{EHSST}, and that $\Ext^*_{A^e}(A,A)$ itself must be
Noetherian.

Although many algebras do not satisfy (Fg.1) and (Fg.2), we can show
that Snashall-Solberg theory works for block algebras of a
transporter category algebra $k\GP$. From \cite{X3} we know that
$\Ext^*_{k\GP}(\k,\k)$ is a Noetherian graded commutative ring such
that $\Ext^*_{k\GP}(\m,\n)$ is finitely generated over it, for any
pair of $\m,\n\in k\GP$-mod. It was also showed there that
$\Ext^*_{k\GP^e}(k\GP,k\GP)$ is Noetherian. We shall prove that
$\Ext^*_{k\GP}(\k,\k)$- and $\Ext^*_{k\GP^e}(k\GP,k\GP)$-actions on
$\Ext^*_{k\GP}(\m,\n)$ are compatible and hence it implies that
$\Ext^*_{k\GP}(\m,\n)$ is finitely generated over the Hochschild
cohomology ring as well.

\begin{theorem} Let $\m\in k\C$-mod. We have a commutative diagram
$$
\xymatrix{\Ext^*_{kF(\C)}(\k,\k) \ar[rr]^{\cong} \ar[d] && \Ext^*_{k\C}(\k,\k) \ar[d]^{-\hotimes\m}\\
\Ext^*_{k\C^e}(k\C,k\C) \ar[rr]_(.55){-\otimes_{k\C}\m} &&
\Ext^*_{k\C}(\m,\m)}
$$
with the left vertical map an injective algebra homomorphism.
\end{theorem}

This actually generalizes \cite[Theorem A]{X1}. Let us first recall
some other necessary constructions from \cite{X1}. For any category
$\C$ there is a \textit{category of factorizations in} $\C$, written
as $F(\C)$. The objects are the morphisms in $\C$. When a morphism
$\alpha\in\Mor\C$ is regarded as an object in $F(\C)$, we will
denote it by $[\alpha]$ to distinguish their roles. If $[\alpha],
[\beta]$ are two objects in $F(\C)$, then a morphism
$[\alpha]\to[\beta]$ is a pair $(\mu,\gamma)$,
$\mu,\gamma\in\Mor\C$, such that $\beta=\mu\alpha\gamma$ (that is,
$\alpha$ is a factor of $\beta$). When $\C$ is a group, $F(\C)$
plays the role of the diagonal subgroup $\triangle G \subset G\times
G$. Indeed, there is a category equivalence $\triangle G\simeq
F(G)$.

Given a morphism $\alpha$ in $\C$, we denote by $t(\alpha)$ and
$s(\alpha)$ the target and source of $\alpha$. They induce two
functors $t : F(\C)\to\C$ and $\nabla=(t,s) :
F(\C)\to\C^e=\C\times\C^{op}$, fitting into the following
commutative diagram
$$
\xymatrix{F(\C) \ar[rr]^{\nabla} \ar[dr]_{t} && \C^e \ar[dl]^{p}\\
&\C& ,}
$$
where $p$ is the projection. By definition $t$ and $s$ send
$[\alpha]$ to the target and source of $\alpha$, respectively. In
\cite{X1}, we investigated the left Kan extensions $LK^{\C}_{F(\C)}
: kF(\C)$-mod $\to k\C$-mod and $LK^{\C^e}_{F(\C)} : kF(\C)$-mod
$\to k\C^e$-mod, proving $LK^{\C}_{F(\C)}\k\cong\k$ and
$LK^{\C^e}_{F(\C)}\k\cong k\C$. Furthermore $LK^{\C}_{F(\C)}$
induces an isomorphism $\Ext^*_{kF(\C)}(\k,\k) \to
\Ext^*_{k\C}(\k,\k)$, while $LK^{\C^e}_{F(\C)}$ induces an injective
algebra homomorphism
$\Ext^*_{kF(\C)}(\k,\k)\to\Ext^*_{k\C^e}(k\C,k\C)$. Especially
$F(\C)$ is connected if and only if $\C$ is. At last $LK^{\C}_{\C^e}
: k\C^e$-mod $\to k\C$-mod is explicitly expressed as
$LK^{\C}_{\C^e}\cong-\otimes_{k\C}\k$, where $\k$ is the trivial
left $k\C$-module.

When $\m=\k$, the lower horizontal map becomes the split surjection
in \cite[Theorem A]{X1}.

\begin{proof} Consider the bar resolution $\B^{F(\C)}_*\to \k\to 0$,
and a map $f : \B^{F(\C)}_n=\c_n(\Id_{F(\C)}/-,k) \to \k$
representing an element $\xi\in\Ext^n_{kF(\C)}(\k,\k)$. We need to
prove
$LK^{\C^e}_{F(\C)}\xi\otimes_{k\C}\m=LK^{\C}_{F(\C)}\xi\hotimes\m
\in\Ext^n_{k\C}(\m,\m)$. We do it by explicit calculations.

Firstly, $LK^{\C^e}_{F(\C)}f :
LK^{\C^e}_{F(\C)}\B^{F(\C)}_n=\c_n(\nabla/-,k)\to
LK^{\C^e}_{F(\C)}\k=k\C$ is given on each $(x,y)\in\Ob\C^e$ as
$$
(LK^{\C^e}_{F(\C)}f)_{(x,y)} : \c_n(\nabla/(x,y),k) \to
k\C(x,y)=k\Hom_{\C}(y,x)
$$
by
$$
([\alpha_*],(\beta_*,\gamma_*)) \mapsto
f_{[\alpha_*]}\beta_n\alpha_n\gamma_n
$$
Here we denote by
$([\alpha_*],(\beta_*,\gamma_*))=([\alpha_0],(\beta_0,\gamma_0))\to
\cdots\to([\alpha_n],(\beta_n,\gamma_n))$ a base element of
$\c_n(\nabla/(x,y),k)$, where $(\beta_i,\gamma_i) :
\nabla([\alpha_i])\to (x,y)$ is a morphism in $\C^e$. From
$([\alpha_*],(\beta_*,\gamma_*))$ we can extract a base element of
$\c_n(\Id_{F(\C)}/[\beta_n\alpha_n\gamma_n],k)$, written as
$[\alpha_*]=[\alpha_0]\to \cdots\to[\alpha_n]$, so our definition
makes sense. Note that
$\beta_0\alpha_0\gamma_0=\cdots=\beta_n\alpha_n\gamma_n$.

Secondly, in a similar fashion, $LK^{\C}_{F(\C)} f :
LK^{\C}_{F(\C)}\B^{F(\C)}_n=\c_n(t/-,k)\to LK^{\C}_{F(\C)}\k =\k$ is
given on each $x\in\Ob\C$ as
$$
(LK^{\C}_{F(\C)} f)_x : \c_n(t/x,k) \to \k(x)=k
$$
by
$$
([\alpha_*],\mu_*) \mapsto f_{[\alpha_*]}.
$$
Here we denote by
$([\alpha_*],\mu_*)=([\alpha_0],\mu_0)\to\cdots\to([\alpha_n],\mu_n)$
a base element of $\c_n(t/x,k)$ and $[\alpha_*]$ the encoded base
element of $\c_n(\Id_{F(\C)}/[\mu_n\alpha_n],k)$.

Thirdly, $LK^{\C^e}_{F(\C)}\xi\otimes_{k\C}\m$ is represented by
$LK^{\C^e}_{F(\C)}f\otimes_{k\C}\Id_{\m}$, while
$LK^{\C}_{F(\C)}\xi\hotimes\m$ is represented by some
$(LK^{\C}_{F(\C)} f\hotimes\Id_{\m})\circ\Phi_n$ provided that $\P_n
\to \m \to 0$ is a projective resolution and $\Phi_* : \P_n \to
LK^{\C}_{F(\C)}\B^{F(\C)}_*\hotimes\m$ is a lifting of the identity
map of $\m$. Let us take $\P_* =
LK^{\C^e}_{F(\C)}\B^{F(\C)}_*\otimes_{k\C}\m$ and construct
$\Phi_*$, which has to be unique up to homotopy. To this end, for
any $x\in\Ob\C$ and $n\ge 0$ we define
$$
\Phi^x_n : 1_x\cdot LK^{\C^e}_{F(\C)}\B^{F(\C)}_n\otimes_{k\C}\m \to
1_x\cdot LK^{\C}_{F(\C)}\B^{F(\C)}_*\otimes\m(x)
$$
by
$$
([\alpha_*],(\beta_*,\gamma_*))\otimes(m_w)\mapsto
([\alpha_*],\beta_*)\otimes(\beta_n\alpha_n\gamma_n)\cdot m_y.
$$
Here we write each element in $\m=\bigoplus_{w\in\Ob\C}\m(w)$ as
$\sum_w m_w$, in which $m_w\in\m(w)$. Directly from the definition
one can verify
$LK^{\C^e}_{F(\C)}f\otimes_{k\C}\Id_{\m}=(LK^{\C}_{F(\C)}
f\hotimes\Id_{\m})\circ\Phi$. It means
$LK^{\C^e}_{F(\C)}\xi\otimes_{k\C}\m=LK^{\C}_{F(\C)}\xi\hotimes\m
\in\Ext^n_{k\C}(\m,\m)$.
\end{proof}

To consider Snashall-Solberg's theory, we want the algebra in
question to be indecomposable. Suppose $k\GP=\prod_i b_i$ is a
decomposition into (indecomposable) block algebras, see for instance
\cite[Section 13]{A}. Since $k\GP$ is Gorenstein, so are its blocks.
Let $b$ be a block of $k\GP$. The above theorem implies that
$\Ext^*_b(\m,\n)$ is a finitely generated module over the Noetherian
ring $\Ext^*_{b^e}(b,b)$, if $\m,\n\in b$-mod. It means
Snashall-Solberg's theory works perfectly for blocks of finite
transporter category algebras.

Here we record some standard properties from Snashall-Solberg's
theory. For convenience, we write
$$
V_b(\m,\n)=V_{\HH(b)}(\m,\n)
$$
if $b$ is a block of $k\GP$ and $\m, \n\in b$-mod. To be consistent,
write $V_b(\m)=V_b(\m,\m)$. Some terminologies are recalled first.

\begin{definition} Let $A$ be a finite-dimensional algebra and $P_* \to M \to 0$
the minimal projective resolution of $M$. Then the
\textit{complexity} of $M$ is
$$
c_A(M) = \min\{s \in \mathbb{N} \bigm{|} r\in \mathbb{R} \
\mbox{such that}\ \dim_k P_n\le r n^{s-1},\ \mbox{for}\ n \gg 0 \}
$$
\end{definition}

Let $(-)^{\wedge}=\Hom_k(-,k) : k\C$-mod $\to k\C^{op}$-mod be the
$k$-dual functor. Recall that $B\C\simeq B\C^{op}$ and thus
$\Ext^*_{k\C}(\k,\k)\cong\Ext^*_{k\C^{op}}(\k,\k)$. It implies that
$\Ext^*_{k\C}(\k,\k)$ also acts on
$\Ext^*_{k\C^{op}}(\n^{\wedge},\m^{\wedge})$, for any $\m,\n\in
k\GP$-mod. Some of these constructions pass to every block algebra
of $k\C$.

The following statements are taken from Snashall-Solberg \cite{SS}
and Erdmann-Holloway-Snashall-Solberg-Taillefer \cite{EHSST}. Note
that a block of a transporter category is Gorenstein, but in general
neither selfinjective nor symmetric.

\begin{theorem} Let $G$ be a finite group, $\P$ a finite $G$-poset, and $k$
an algebraically closed field of characteristic $p$ dividing the
order of $G$. Suppose $b$ is a block of $k\GP$ and $\m,\n$ are two
finitely generated modules of $b$. Then
\begin{enumerate}
\item $V_b(\m)=\bigcup_{\mathfrak{S}}V_b(\m,\mathfrak{S})=\bigcup_{\mathfrak{S}}V_b(\mathfrak{S},\m)$,
where $\mathfrak{S}$ runs over the set of all simple $b$-modules.

\item If $0\to \m_1 \to \m_2\to\m_3\to 0$ is an exact sequence, then
$V_b(\m_i)\subset V_b(\m_j)\cup V_b(\m_l)$ for
$\{i,j,l\}=\{1,2,3\}$.

\item $V_b(\m_1\oplus\m_2)=V_b(\m_1)\cup V_b(\m_2)$.

\item $V_b(\m)=V_b(\Omega^n(\m))$ for any integer $n$ such
that $\Omega^n(\m)\ne 0$.

\item $V_b(\m,\n)=V_{b^{op}}(\n^{\wedge},\m^{\wedge})$.
Particularly $V_b(\m)=V_{b^{op}}(\m^{\wedge})$.

\item $\dim V_b(\m)=c_b(\m)$.

\item $V_b(\m)$ is trivial if and only if $\m$ has
finite projective dimension if and only if $\m$ has finite injective
dimension.

\item Let $\mathfrak{a}$ be a homogeneous ideal in
$\HH(b)=\Ext_{b^e}(b,b)$. Then there exists a module
$\m_{\mathfrak{a}}\in b$-mod such that
$V_b(\m_{\mathfrak{a}})=V(\mathfrak{a})$.

\item If $V_b(\m)\cap V_b(\n)$ is trivial, then
$\Ext^i_b(\m,\n)=0$ for $i > {\rm inj.dim}b$.

\item If $V_b(\m)=V_1\cup V_2$ for some homogeneous non-trivial
varieties $V_1$ and $V_2$ with $V_1\cap V_2$ trivial, then
$\m=\m_1\oplus\m_2$ with $V_b(\m_1)=V_1$ and $V_b(\m_2)=V_2$.
\end{enumerate}

\begin{proof} The first five properties follow directly from the
definition of support variety for a module of Snashall-Solberg, and
can be found in \cite{SS}.

The rest come from \cite[Theorem 2.5, Theorem 4.4, Proposition 7.2,
Theorem 7.3]{EHSST}.
\end{proof}
\end{theorem}

\section{Standard properties of $V_{\gp}(\m)$}

Snashall-Solberg's theory on a block algebra of a transporter
category algebra is quite satisfactory in many ways. However,
Hochschild cohomology rings do not behave well comparing with
ordinary cohomology rings. For example, since an algebra
homomorphism does not necessarily induce a homomorphism between
their Hochschild cohomology rings, certain important properties in
Carlson's theory are not expected to exist in Snashall-Solberg's
theory. This is one of the reasons why we believe $V_{\gp}(\m)$ has
various advantages over $V_b(\m)$ which we try to demonstrate in the
rest of this paper.

\subsection{Principal block}

Let $k\GP$ be a transporter category algebra. Remember that we
assume $\gp$ is connected, which is equivalent to saying that $\k$
is indecomposable. We pay special attention to a special block of
the transporter category algebra, closely related to our support
variety theory.

\begin{definition} Assume $\C$ is a finite connected category.
We call a block of $k\C$ the \textit{principal block} if the
(indecomposable) trivial module $\k$ belongs to it, and consequently
name the block $b_0$.
\end{definition}

Since one can take a minimal projective resolution of $\k$
consisting of projective modules belonging to $b_0$,
$\Ext^*_{k\C}(\k,\k)=\Ext^*_{b_0}(\k,\k)$ is an invariant of the
principal block, comparable to the group case, see for example
\cite{L1}.

Return to transporter category algebras. We claim

\begin{enumerate}
\item[(i)] $\Ext^*_{k\GP}(\k,\k)$ is a (usually proper) subring of
$\Ext^*_{b_0^e}(b_0,b_0)$;

\item[(ii)] Snashall-Solberg's theory is valid for the subring
$$
\mathbb{H}:=\langle Z(b_0),\Ext_{k\GP}(\k,\k)\rangle
\hookrightarrow\Ext_{b_0^e}(b_0,b_0);
$$

\item[(iii)] $\Ext^*_{k\GP}(\k,\k)\hookrightarrow\langle
Z(b_0),\Ext^*_{k\GP}(\k,\k)\rangle$ induces an isomorphism after
quotient out nilpotent elements.
\end{enumerate}

The ordinary cohomology ring is known to be finitely generated
\cite{Q1,X2}, and thus so are the rings $\langle
Z(b_0),\Ext^*_{k\GP}(\k,\k)\rangle$ and $\mathbb{H}$. Assuming these
claims, along with Theorem 3.4.2, it guarantees that
Snashall-Solberg's theory can be implemented to
$\mathbb{H}\subset\HH(b_0)$ for the principal block of $k\GP$.

Claim (ii) comes from Theorem 3.4.2. Claim (iii) is easy to verify
since the commutative local algebra $Z(b_0)$ quotients out the
unique maximal ideal, that is, its nilradical, is exactly $k$ for
$k$ is algebraically closed. To establish Claim (i), we have the
following proposition.

\begin{proposition} Let $\C$ be a finite connected category and $b_0$
the principal block of $k\C$. The injective homomorphism
$$
\Ext^*_{k\C}(\k,\k)\hookrightarrow\Ext^*_{k\C^e}(k\C,k\C)
$$
induces an injective homomorphism
$$
\Ext^*_{k\C}(\k,\k)\hookrightarrow\Ext^*_{b_0^e}(b_0,b_0).
$$

\begin{proof} As a $k\C^e$-module, $k\C=\oplus_i b_i$. Since $LK^{\C}_{\C^e} k\C\cong\k$ and $LK^{\C}_{\C^e}$
preserves direct sums, we see $LK^{\C}_{\C^e}
b_0=b_0\otimes_{k\C}\k=\k$, and $LK^{\C}_{\C^e} b_i=0$ if $b_i$ is
not principal. It implies that the split surjection in \cite{X1}
$$
\Ext^*_{k\C^e}(k\C,k\C)\to\Ext^*_{k\C}(\k,\k)
$$
induced by $LK^{\C}_{\C^e}$, restricts to a split surjection
$\Ext^*_{b_0^e}(b_0,b_0)\to\Ext^*_{k\C}(\k,\k)$. Hence we are done.
\end{proof}
\end{proposition}

From the preceding discussions, we get
$$
V_{\gp}(\m)=V_{\mathbb{H}}(\m)
$$
for all $\m\in b_0$-mod. Note that, for any $\m\in b_0$-mod, there
exists a finite surjective map
$$
V_{b_0}(\m) \to V_{\gp}(\m).
$$
We do not know yet when it becomes an isomorphism.

\subsection{Standard properties}

A subvariety of $V_{\gp}$ is \textit{trivial} if it is the ideal
consisting of all positive degree elements of $\Ext_{k\GP}(\k,\k)$.

\begin{theorem} Let $G$ be a finite group, $\P$ a finite $G$-poset, and $k$
an algebraically closed field of characteristic $p$ dividing the
order of $G$. Suppose $\m,\n$ are two finitely generated modules of
$k\GP$. Then
\begin{enumerate}
\item $V_{\gp}(\m)=\bigcup_{\mathfrak{S}}V_{\gp}(\m,\mathfrak{S})=
\bigcup_{\mathfrak{S}}V_{\gp}(\mathfrak{S},\m)$, where
$\mathfrak{S}$ runs over the set of all simple $k\GP$-modules.

\item If $0\to \m_1 \to \m_2\to\m_3\to 0$ is an exact sequence, then
$V_{\gp}(\m_i)\subset V_{\gp}(\m_j)\cup V_{\gp}(\m_l)$ for
$\{i,j,l\}=\{1,2,3\}$.

\item $V_{\gp}(\m_1\oplus\m_2)=V_{\gp}(\m_1)\cup V_{\gp}(\m_2)$.

\item $V_{\gp}(\m)=V_{\gp}(\Omega^n(\m))$ for any integer $n$ such
that $\Omega^n(\m)\ne 0$.

\item $V_{\gp}(\m,\n)=V_{{\GP}^{op}}(\n^{\wedge},\m^{\wedge})$.
Particularly $V_{\gp}(\m)=V_{{\GP}^{op}}(\m^{\wedge})$.

\item $\dim V_{\gp}(\m)=c_{k\GP}(\m)$.

\item $V_{\gp}(\m)$ is trivial if and only if $\m$ has
finite projective dimension if and only if $\m$ has finite injective
dimension.

\item Let $\mathfrak{a}$ be a homogeneous ideal in
$\Ext_{k\GP}(\k,\k)$. Then there exists a module
$\m_{\mathfrak{a}}\in k\GP$-mod such that
$V_{\gp}(\m_{\mathfrak{a}})=V(\mathfrak{a})$.

\item[(9a)] If $V_{\gp}(\kappa_{\Hom_k(M,N)})$ is trivial, then
$\Ext^i_{k\GP}(\kappa_M,\kappa_N)=0$ for $i > \dim\P$.

\item[(9b)] If $\m, \n \in b_0$-mod and $V_{\gp}(\m)\cap V_{\gp}(\n)$ is trivial, then
$\Ext^i_{k\GP}(\m,\n)=0$ for $i > {\rm inj.dim}b_0$.

\item[(10a)] If $\chi(\P) \equiv 1$ (mod $p$) and $V_{\gp}(\kappa_M)=V_1\cup V_2$
for some homogeneous non-trivial
varieties $V_1$ and $V_2$ with $V_1\cap V_2$ trivial, then
$\kappa_M=\kappa_{M_1}\oplus\kappa_{M_2}$ with
$V_{\gp}(\kappa_{M_1})=V_1$ and $V_{\gp}(\kappa_{M_2})=V_2$.

\item[(10b)] If $\m \in b_0$-mod such that $V_{\gp}(\m)=V_1\cup V_2$
for some homogeneous non-trivial
varieties $V_1$ and $V_2$ with $V_1\cap V_2$ trivial, then
$\m=\m_1\oplus\m_2$ with $V_{\gp}(\m_1)=V_1$ and
$V_{\gp}(\m_2)=V_2$.
\end{enumerate}

\begin{proof} The first five properties follow directly from the
definition of $V_{\gp}(\m,\n)$. The proofs are entirely analogous to
the group case, see \cite{B2, Evens}.

The proof of (6) is exactly the same as that for groups, see for
example \cite[Proposition 5.7.2]{B2}. As for (7), one direction is
straightforward. Also by Lemma 2.3.1 (2), it is equivalent to saying
that $\m$ is of finite injective dimension. Now let us assume
$V_{\gp}(\m)$ is trivial. Then by (6) $c_{k\GP}(\m)=0$. It forces
the minimal projective resolution of $\m$ to be finite.

Part (8) follows from $V_{\gp}(\m)=V_{\mathbb{H}}(\m)$ if $\m\in
b_0$-mod. Let $\mathfrak{a}$ be a homogeneous ideal in
$\Ext_{k\GP}(\k,\k)$. Then we can define a homogeneous ideal
$\mathfrak{a}'=\langle\Rad{Z(b_0)},\mathfrak{a}\rangle$ of
$\mathbb{H}$. From Snashall-Solberg's theory \cite[Theorem
4.4]{EHSST}, the general form of Theorem 3.4.4 (8) for $b_0$, there
exists a $b_0$-module $\m_{\mathfrak{a}}$ such that
$V_{\mathbb{H}}(\m_{\mathfrak{a}})=V(\mathfrak{a}')$. But
$V(\mathfrak{a}')$ is identified with $V(\mathfrak{a})$ under the
isomorphism $V_{\mathbb{H}} \to V_{\gp}$, and
$V_{\gp}(\m_{\mathfrak{a}})=V_{\mathbb{H}}(\m_{\mathfrak{a}})$.

Since (9b) and (10b) are Theorem 3.4.4 (9) and (10) specialized to
the principal block, we shall prove only (9a) and (10a).

To prove (9a), we notice that
$\Ext^i_{k\GP}(\kappa_M,\kappa_N)\cong\Ext^i_{k\GP}(\k,\kappa_{\Hom_k(M,N)})$
(see Section 2.5). From (7), the assumption implies that
$\kappa_{\Hom_k(M,N)}$ has finite injective dimension. Hence the
statement follows from Lemma 2.3.1.

As to (10a), we recall from Section 3.2 that if $\chi(\P) \equiv 1$
(mod $p$) then there exists a finite surjective map
$V_{\gp}(\kappa_M) \to V_G(M)$ for $M \in kG$-mod. Under the
circumstance, let $V'_i$ be the images of $V_i$, for $i = 1, 2$, we
see that $V_G(M)=V'_1\cup V'_2$ with $V'_1\cap V'_2$ trivial. By
\cite[Theorem 5.12.1]{B2}, $M=M_1\oplus M_2$ satisfying
$V_G(M_1)=V'_1$ and $V_G(M_2)=V'_2$. Then $\kappa_M=
\kappa_{M_1}\oplus\kappa_{M_2}$ and hence
$V_{\gp}(\kappa_M)=V_{\gp}(\kappa_{M_1})\cup V_{\gp}(\kappa_{M_2})$.
Moreover since the preceding map between varieties restricts to
finite surjective maps $V_{\gp}(\kappa_{M_i}) \to V_G(M_i)$ for $i =
1, 2$. It implies $V_i=V_{\gp}(\kappa_{M_i})$ for $i= 1, 2$.
\end{proof}
\end{theorem}

By Theorem 4.2.1 (7) and Lemma 2.3.1, $V_{\gp}(\kappa_M)$ is trivial
if and only if $\kappa_M$ has finite projective dimension if and
only if $M$ is a projective $kG_x$-module for all $x\in\Ob\GP$. If
$\P=\E_p$, the poset of all elementary abelian $p$-subgroups of $G$,
by Chouinard's theorem \cite[Theorem 5.2.4]{B2} it is equivalent to
saying that $M$ is a projective $kG$-module.

It is not known yet whether in (9a) and (10a) we may replace the
constant modules with arbitrary modules.

\subsection{Consequences of module filtrations}

Recall that since $\gp$ is a finite EI-category, every $k\GP$-module
$\m$ is constructed from atomic modules $\m_x$, where
$\m_x(y)\cong\m(x)$ if $y\cong x$ or zero otherwise. The following
result says that only $\m$-objects contribute to the variety
$V_{\gp}(\m)$. Moreover the non-isomorphisms do not play a big role.

\begin{proposition} We have $V_{\gp}(\m)= \bigcup_{[x]}V_{\gp}(\m_x)
=\bigcup_{[x]}\iota^{\gp}_{G\propto[x]}V_{G\propto[x]}(\m_x)$, where
$[x]$ runs over the set of isomorphism classes of $\m$-objects.

Particularly $V_{\gp}(\kappa_M)= \bigcup_{[x]}\iota^{\gp}_{G\propto
[x]}V_{G\propto[x]}(\kappa_M)$, for any $M\in kG$-mod.

\begin{proof} The containment $V_{\gp}(\m)\subset \bigcup_{[x]\in\Is\GP}V_{\gp}(\m_x)$
follows from Theorem 4.2.1 (2). We establish the equality. Firstly
we note that
$\Ext^*_{k\GP}(\m_x,\m_x)\cong\Ext^*_{k(G\propto[x])}(\m_x,\m_x)$
naturally. Moreover we have a commutative diagram
$$
\xymatrix{\Ext^*_{k\GP}(\k,\k) \ar[rr]^{-\hotimes\m}
\ar[d]_{\res^{G\propto\P}_{G\propto[x]}}
&& \Ext^*_{k\GP}(\m,\m) \ar[d]^{\res^{G\propto\P}_{G\propto[x]}}\\
\Ext^*_{k(G\propto[x])}(\k,\k) \ar[rr]_{-\hotimes\m_x} &&
\Ext^*_{k(G\propto[x])}(\m_x,\m_x)}
$$
To establish the above commutative diagram, we can represent
cohomology classes by $n$-fold extensions and notice that both
restrictions and $-\hotimes-$ are exact.

The $\Ext^*_{k\GP}(\k,\k)$-action on $\Ext^*_{k\GP}(\m_x,\m_x)$
factors through the action by $\Ext^*_{G\propto[x]}(\k,\k)$. Hence
we have
$V_{\gp}(\m_x)=\iota^{\gp}_{G\propto[x]}V_{G\propto[x]}(\m_x)$.
Based on the same diagram we see that $I_{\gp}(\m)$ kills
$\Ext^*_{k\GP}(\m_x,\m_x)$. It means $V_{\gp}(\m_x)\subset
V_{\gp}(\m)$.
\end{proof}
\end{proposition}

We remind the reader that since $G\propto[x]$ as a category is
equivalent to $G_x$, $V_{G\propto[x]}(\m_x)$ can be identified with
$V_{G_x\times x}(\m(x))$.

\begin{corollary} We have $V_{\gp}(\m)=\bigcup_{[x]}\iota^{\gp}_{G_x\times x} V_{G_x\times x}(\m(x))$, where
$[x]$ runs over the set of isomorphism classes of $\m$-objects, and
hence $V_{\gp}(\kappa_M)=\bigcup_{[x]}\iota^{\gp}_{G_x\times
x}V_{G_x\times x}(M)$ for any $M\in kG$-mod.

If $G$ acts trivially on a connected poset $\P$, then
$V_{G\times\P}(\m)=\bigcup_x\iota^{\gp}_{G_x\times x}V_G(\m(x))$.
Especially $V_{G\times\P}(\kappa_M)=V_G(M)$ for any $M\in kG$-mod.
\end{corollary}

The reader may go back and have another look at Example 3.1.2. We
note that if $x\cong y$ in $\gp$, then there exists an element $g\in
G$ inducing an isomorphism by conjugation $G_x \to G_y$. It implies
that $V_{G_x\times x}(\m(x))\cong V_{G_y\times y}(\m(y))\cong
V_{G\propto [x]}(\m_x)$ and
$$
\iota^{\gp}_{G_x\times x}V_{G_x\times
x}(\m(x))=\iota^{\gp}_{G_y\times y}V_{G_y\times
y}(\m(y))=\iota^{\gp}_{G\propto[x]}V_{G\propto[x]}(\m_x)=V_{\gp}(\m_x).
$$

\begin{corollary} If $\m\in k\GP$-mod and $H\subset G$ is a
subgroup, then $\iota^{\gp}_{H\propto\P}V_{H\propto\P}(\m)\subset
V_{\gp}(\m)$.

\begin{proof} By the preceding result, we have $V_{H\propto\P}(\m)=
\bigcup_x\iota^{H\propto\P}_{H_x\times x} V_{H_x\times x}(\m(x))$,
and $V_{\gp}(\m)=\bigcup_x\iota^{\gp}_{G_x\times x} V_{G_x\times
x}(\m(x))$, where $x$ runs over the set of all $\m$-objects.

Since $\iota^{\gp}_{H\propto\P}\iota^{H\propto\P}_{H_x\times
x}=\iota^{\gp}_{G_x\times x}\iota^{G_x\times x}_{H_x\times x}$ and
$\iota^{G_x\times x}_{H_x\times x}V_{H_x\times x}(\m(x))\subset
V_{G_x\times x}(\m(x))$ is known in Carlson's theory, our claim
follows.
\end{proof}
\end{corollary}

In a similar fashion we can show if $\Q$ is a $G$-subposet of $\P$
then
$$
\iota^{\gp}_{G\propto\Q}V_{H\propto\Q}(\m)\subset V_{\gp}(\m)
$$
for any $\m\in k\GP$-mod.

\subsection{Varieties of Mackey functors}

There are various ways to construct modules for a transporter
category algebra. In \cite{X3} we examined a context of three
categories
$$
\xymatrix{&\gp \ar[dl]_{\pi} \ar[dr]^{\rho}&\\
G && \C,}
$$
where $\C$ is a quotient category of $\gp$. An illuminating example
is that, when $\P$ is a poset of subgroups and $G$ acts via
conjugation, we may define an \textit{orbit category}
$\C=\O_{\P}(G)$ as the quotient category of $\gp$ such that
$\Hom_{\O_{\P}(G)}(H,K)=K\backslash\Hom_{\gp}(H,K)=K\backslash
N_G(H,K)$. Thus every $k\O_{\P}(G)$-module (resp.
$k\O^{op}_{\P}(G)$)-module) can be restricted to a $k\GP$-module
(resp. $k\GP^{op}$-module). If $\S$ is the poset of all subgroups of
$G$, a $k\O^{op}_{\S}(G)$-module is often called a
\textit{coefficient system} (for equivariant cohomology of
$G$-spaces). The natural examples are the \textit{Mackey functors}
\cite{tD, T, W1}. A Mackey functor may be defined as a
$k\O_{\S}(G)$-bimodule enjoying certain extra properties. However
after restriction (brutal truncation) we have no trouble to get a
$k\O_{\P}(G)$-bimodule for an arbitrary $G$-subposet $\P\subset\S$.
Thus every Mackey functor gives rise to a $k\GP$-bimodule. In light
of this it is tempting to introduce the following concept.

\begin{definition} Let $\m$ be a Mackey functor and $\P$ a $G$-poset
consisting of some subgroups of $G$. The variety of the Mackey
functor is $V_{\gp}(\m)$, where $\m$ is regarded as a (either left
or right) module over $k\GP$.
\end{definition}

The reason why choosing the left or right $k\GP$-module structure
for a Mackey functor shall not make a difference is that, as far as
we concern, only its values on objects matter, see Proposition
4.3.1.

\begin{example}
Let us now consider the poset $\E_p$ and the Mackey functor
$\mathfrak{H}^n_M:=\H^n(-,M)$, where $n$ is a non-negative integer
and $M$ is a $kG$-module. Then we have
$$
V_{G\propto\E_p}(\mathfrak{H}^n_M)=\cup_E
\iota^{G\propto\E_p}_{N_G(E)}V_{N_G(E)}(\H^n(E,M)).
$$
Particularly when $n=0$ and $M=k$, $\mathfrak{H}^0_k=\k$. Based on
discussions in Section 3.2 we get
$$
V_{G\propto\E_p}(\mathfrak{H}^0_k)=V_G.
$$
\end{example}

It shall be interesting to understand the varieties of various
Mackey functors.

\section{Further properties}

In this section we shall deal with comparing varieties of
categories. The main results are the generalized Quillen
stratification and its consequences.

Bear in mind that, for the sake of simplicity, if $\D\to\C$ is a
functor and $\m$ is a $k\C$-module, when it will not cause
confusions, we often regard $\m$ (instead of writing
$\Res^{\C}_{\D}\m$) as a $k\D$-module, even though its underlying
vector space structure usually changes.

\subsection{Quillen stratification}

In \cite{Q1, Q2} Quillen worked with $G$-spaces and equivariant
cohomology rings. Here we are interested in $G$-spaces which are
classifying spaces of finite $G$-posets. In order to make a smooth
transition from $G$-spaces and equivariant cohomology to $G$-posets
and transporter category cohomology, we first recall some of the
original constructions and then restrict to our case.

Given a $G$-space $X$ and a $H$-space $Y$, a morphism from $(H,Y)\to
(G,X)$ is a pair $(\phi, F)$ such that $\phi : H\to G$ is a group
homomorphism and $\theta : Y \to X$ is a continuous map satisfying
the condition that
$$
\phi(h)\theta(y)=\theta(hy),\ \ \forall h\in H, y\in Y.
$$
Such a morphism induces a continuous map $\Theta : EH\times_H Y\to
EG\times_G X$ and thus a restriction map between equivariant
cohomology
$$
\res^{G,X}_{H,Y} : \H^*_{G}(X) \to \H^*_H(Y).
$$
One can compare these constructions with those introduced before
Example 3.1.2.

In Quillen's papers \cite{Q1,Q2}, it is proved that if $X$ is a
$G$-space, either compact or paracompact with finite cohomological
dimension, then there exists an F-isomorphism
$$
\q=\q_{G,X} : \H^*_G(X)\to{\lim}_{\A_p(X)}\H^*(E)
$$
where $\A_p(X)$ is called the \textit{Quillen category} and its
objects will be named \textit{Quillen pairs} for the $G$-space $X$
in this article. More explicitly, the objects in $\A_p(X)$ are of
the form $(E,C)$, where $E$ is an elementary abelian $p$-subgroup of
$G$ and $C$ is a (non-empty) connected component of $X^E$. A
morphism from $(E',C') \to (E,C)$ is a pair $(\phi, \i)$ with
$\phi=c_g : E' \to E$ for some element $g\in G$ and $\i : gC \to C'$
an inclusion. Particularly if $g$ can be chosen to be the identity
element of $G$, then we call $(E',C')$ a \textit{Quillen subpair} of
$(E,C)$, and write $(E',C')\le (E,C)$. In this way, all Quillen
pairs form a poset $\Q_p$, with obvious $G$-action.

\begin{remark}
Another way to construct the Quillen category is to define it as a
quotient category of the transporter category $G\propto\Q_p$, where
$\Q$ is the $G$-poset of all Quillen pairs for $X$, through
$$
\Hom_{\A_p(X)}((E',C'),(E,C))=\Hom_{G\propto\Q_p}((E',C'),(E,C))/C_G(E',C'),
$$
where $C_G(E',C')=\{g\in G \bigm{|} gC'=C', ghg^{-1}=h\ \mbox{for
all}\ h \in E\}$ is called the \textit{centralizer} of $(E',C')$. We
also denote by $N_G(E',C')=G_{(E',C')}$ the normalizer of $(E',C')$
and $W_G(E',C')=N_G(E',C')/C_G(E',C')$ the \textit{Weyl group} of
$(E',C')$.
\end{remark}

Quillen's map $\q$ is induced by $(E,\bullet) \to (E, C) \to (G,X)$,
where $\bullet$ is sent into $C$. Since $E$ acts trivially on $C$,
the choice of an embedding $\bullet\to C$ does not matter. In fact
in any case this map induces a surjective map
$\H^*_E(C)\to\H^*_E(\bullet)$ which becomes an isomorphism after
quotient out nilpotents in both rings. Consequently
$\Spec\H_E(C)=\Spec\H_E(\bullet)$, that is,
$V_{E,C}=V_{E,\bullet}=V_E$.

Let $V_{G,X}=\Spec\H_G(X)$. The geometric version of Quillen's map
is
$$
V_{G,X}=\bigcup_{(E,C)}\iota^{G,X}_{E,C}
V_{E,C}=\bigcup_{(E,C)}\iota^{G,X}_{E,C}\iota^{E,C}_{E,\bullet}
V_{E,\bullet},
$$
where $\iota^{G,X}_{E,C}: V_{E,C}\to V_{G,X}$ is induced by
$\res^{G,X}_{E,C}$. Based on the observation that each morphism
$(E',C')\to (E,C)$ in the category $\A_p(X)$ induces a ring
homomorphism $\H^*_E \to \H^*_{E'}$, which is compatible with the
two maps $\H^*_G(X)\to\H^*_{E'}$ and $\H^*_G(X)\to\H^*_E$ coming
from $(E',C')\to(G,X)$ and $(E,C)\to (G,X)$, Quillen then continued
to demonstrate that $V_{G,X}$ is a disjoint union of some locally
closed subvarieties ${V^{G,X}_{E,C}}^+$ by examining more closely
the following diagram
$$
\xymatrix{V_{E'}=V_{E',C'} \ar[rr]^{\iota^E_{E'}}
\ar[dr]_{\iota^{G,X}_{E',C'}}&&
V_{E,C}=V_E \ar[dl]^{\iota^{G,X}_{E,C}}\\
&V_{G,X}&.}
$$
Here the horizontal map is $\iota^E_{E'}$ (corresponding to $\H^*_E
\to \H^*_{E'}$), \textit{not} the senseless $\iota^{E,C}_{E',C'}$,
since a morphism in $\A_p(X)$ is different from the morphisms
introduced in the second paragraph of this section.

\begin{theorem}[Quillen] The dimension of $V_{G,X}$ equals the maximum
$p$-rank of an elementary abelian $p$-subgroup from Quillen pairs.
\end{theorem}

We shall come back to Quillen's results shortly after we sort out
all terminologies.

Now we turn to the case of $X=B\P$ for a finite $G$-poset $\P$. Let
us remind the reader of several relevant results.

\begin{enumerate}
\item[(a)] Suppose a $G$-space $X$ is either compact or paracompact with
finite cohomological dimension. If $\H^*(X)$ is finite-dimensional
then so is $\H^*(C)$ for all possible $C$ from Quillen pairs,
\cite[Corollary 4.3]{Q1}.

\item[(b)] $(B\P)^g=B(\P^g)$ (fixed points) for any $g\in G$, \cite[Proposition 66.3]{CR};

\item[(c)] Let $\S_p$ and $\E_p$ be the $G$-posets of non-identity $p$-subgroups
and of elementary abelian $p$-subgroups of $G$. Then obviously any
elementary abelian $p$-subgroup $E$ has the property that both
$(B\S_p)^E$ and $(B\E_p)^E$ are non-empty.

When $X=B\P$ for a $G$-poset $\P$ and $M\in kG$-mod, we know
$(\Id_G,pt) : (G,B\P) \to (G,\bullet)$ induces the following
restriction map
$$
\res : \Ext^*_{kG}(k,k)=\H^*(G) \to
\H^*_G(B\P)=\Ext^*_{k\GP}(\k,\k),
$$
which is injective if the Euler characteristic $\chi(\P)=\chi(B\P)$
is invertible in the base field, by invoking the Becker-Gottlieb
transfer \cite{X4}. A typical example of such $G$-posets is $\S_p$,
the poset of non-identity $p$-subgroups in $G$, since Brown's
theorem \cite[Corollary 6.7.4]{B2} says that $\chi(\S_p)\equiv 1$
(mod $p$). There are several well known subposets that are
$G$-homotopy equivalent to it, and thus possess the same property.

Consider the canonical map $(\Id_E, pt) : (E,X)\to (E,\bullet)$,
where $E$ is an elementary abelian $p$-group. A more involved result
says that the restriction $\res : \H^*_E(\bullet)\to\H^*_E(X)$ being
injective is equivalent to saying that $X^E$ is non-empty,
\cite[IV.1, Corollary 1]{H}. If $X=B\P$ such that $\chi(\P)\equiv 1$
(mod $p$), by \cite[Capter III, Theorem 4.3]{Bre} we get
$\chi(\P^E)\equiv 1$ (mod $p$). This implies that
$B\P^E\ne\emptyset$ and that $\res$ is injective, matching the
observation we made using the Becker-Gottlieb transfer.
\end{enumerate}

From now on we focus on the case $X=B\P$ for a finite $G$-poset
$\P$. By doing so, we restrict to a special case of Quillen in order
to get rid of topology and unveil the underlying
algebraic/categorical constructions. Most importantly we gain the
freedom to work with varieties of modules.

Let $E$ be an elementary abelian $p$-subgroup of $G$. Because
$B\P^E=B(\P^E)$, a Quillen pair for the $G$-space $B\P$ is of the
form $(E,B\C)$ in which $\C$ is a connected component of the poset
$\P^E$. Thus under the circumstance, we may write each Quillen pair
as $(E,\C)$ with $\C$ a connected subposet of $\P^E$, and call
$(E,\C)$ a \textit{Quillen pair} for the $G$-poset $\P$.

Quillen stratification in our setting says (to be consistent with
terminologies in Section 3.1 we choose to write
$\iota^{\gp}_{E\times\C}$ for $\iota^{G,B\P}_{E,B\C}$ and
accordingly $V_{\gp}$ for $V_{G,B\P}$ etc.)
$$
V_{\gp}=\bigcup_{(E,\C)}\iota^{\gp}_{E\times\C}V_{E\times\C}
=\bigcup_{(E,\C)}\iota^{\gp}_{E\times\C}\iota^{E\times\C}_{E\times\bullet}V_{E\times\bullet}.
$$
Note that $E\propto\C=E\times\C$ is a subcategory of $G\propto\P$,
and $V_{E\times\C}\cong V_{E\times\bullet}$.

\begin{remark} The theorem of
Alperin-Avrunin-Evens \cite[Theorem 8.3.1]{Evens} says that
$$
V_G(M)=\bigcup_{E\le G}\iota^G_E V_E(M),
$$
where $E$ runs over the set of all elementary abelian subgroups of
$G$ and $M$ is a $kG$-module. Since $G_x$ is the isotropy group of
$x\in\Ob\P$, we know $x\in\Ob\P^E$ if and only if $E\subset G_x$.
Thus if $(E,\C)$ is a Quillen pair, any object $x\in\Ob\C$ satisfies
the condition that $E\subset G_x$. Thus
$$
\iota^{G_x\times x}_{E\times x} V_{E\times x} \subset V_{G_x\times
x}=V_{G\propto[x]}
$$
Based on the Alperin-Avrunin-Evens' theorem, we can rewrite
Proposition 4.3.1 and Corollary 4.3.2.
$$
\begin{array}{ll}
V_{\gp}(\m)&=\bigcup_{[x]}\iota^{\gp}_{G\propto[x]}V_{G\propto[x]}(\m_x)\\
&=\bigcup_{y}\iota^{\gp}_{G_y\times y}V_{G_y\times
y}(\m(y))\\
&=\bigcup_{y,E\subset G_y}\iota^{\gp}_{E\times y}V_{E\times
y}(\m(y))\\
&=\bigcup_{(E,\C)}\iota^{\gp}_{E\times\C}V_{E\times\C}(\m).
\end{array}
$$
These shall be useful for us in the sequel. We note that if
$(E,\C)\cong(^gE,g\C)$ for some $g\in G$, that is, two Quillen pairs
are isomorphic in $\gp$, then
$\iota^{\gp}_{E\times\C}V_{E\times\C}(\m)= \iota^{\gp}_{^gE\times
g\C}V_{^gE\times g\C}(\m)$.
\end{remark}

Let $(E,\C)$ be a Quillen pair and $E'$ a subgroup of $E$. Then we
denote by $\C|E'$ the connected component of $\P^{E'}$ which
contains $\C$, determining a unique Quillen subpair
$(E',\C|E')\le (E,\C)$.  One can easily verify
$$
\iota^{E}_{E'}\iota^{E'\times(\C|E')}_{E'\times
x}=\iota^{E\times\C}_{E\times x}\iota^{E\times x}_{E'\times x} :
V_{E'}=V_{E'\times x} \to V_{E\times x}=V_E,
$$
for every $x\in\Ob\C$.

\begin{lemma} For any
Quillen pair $(E,\C)$ and $\m\in k(E\times\C)$-mod
$$
\iota^{E}_{E'}V_{E'\times(\C|E')}(\m)=\{\iota^{E}_{E'}
V_{E'\times(\C|E')}\}\cap V_{E\times\C}(\m).
$$

\begin{proof} By Corollary
4.3.2,
$V_{E\times\C}(\m)=\bigcup_{x\in\Ob\C}\iota^{E\times\C}_{E\times
x}V_{E\times x}(\m(x))$ because the isotropy subgroup of every
$x\in\C$ is exactly $E$. Then
$$
\begin{array}{ll}
\{\iota^{E}_{E'}V_{E'\times(\C|E')}\}\cap \iota^{E\times\C}_{E\times
x}V_{E\times x}(\m(x))&=\iota^{E\times\C}_{E\times
x}\{[\iota^{E\times x}_{E'\times x}V_{E'\times x}]\cap V_{E\times x}(\m(x))\}\\
&=\iota^{E\times\C}_{E\times
x}\iota^{E\times x}_{E'\times x}V_{E'\times x}(\m(x))\\
&=\iota^{E}_{E'}\iota^{E'\times(\C|E')}_{E'\times x}V_{E'\times
x}(\m(x))
\end{array}
$$
for every $x\in\Ob\C.$ The second equality uses \cite[Proposition
5.7.7]{B2}. Consequently we get
$$
\begin{array}{ll}
\{\iota^{E}_{E'} V_{E'\times(\C|E')}\}\cap
V_{E\times\C}(\m)&=\{\iota^{E}_{E'}
V_{E'\times(\C|E')}\}\cap\{\bigcup_{x\in\Ob\C}
\iota^{E\times\C}_{E\times
x}V_{E\times x}(\m(x))\}\\
&=\bigcup_{x\in\Ob\C}\{\iota^{E}_{E'}\iota^{E'\times(\C|E')}_{E'\times
x}V_{E'\times x}(\m(x))\}\\
&=\iota^{E}_{E'}V_{E'\times(\C|E')}(\m).
\end{array}
$$
The last equality is true since $\m(y)=0$ if $y\in\Ob(\C|E')-\Ob\C$,
by definition.
\end{proof}
\end{lemma}

We set
$$
V^+_{E\times\C}=V_{E\times\C}-\bigcup_{E'\subsetneqq
E}\iota^{E}_{E'}V_{E'\times(\C|E')}
$$

$$
V_{E\times\C}(\m)^+=V^+_{E\times\C}\cap V_{E\times\C}(\m)
$$

$$
{V^{\gp}_{E\times\C}}^+=\iota^{\gp}_{E\times\C}V_{E\times\C}-\bigcup_{E'\subsetneqq
E}\iota^{\gp}_{E'\times(\C|E')}V_{E'\times(\C|E')}
$$

$$
V^{\gp}_{E\times\C}(\m)^+=\iota^{\gp}_{E\times\C}V_{E\times\C}(\m)^+
$$

The following is established by Quillen for arbitrary $\P$ and
$\m=\k$ \cite{Q2} and by Avrunin-Scott for $\P=\bullet$  and
arbitrary modules \cite{AS}. Our generalization is based on both of
these special cases.

\begin{theorem}[Stratification] Let $\gp$ be a finite transporter category,
$k$ an algebraically closed field of characteristic $p$ dividing the
order of $G$ and $\m\in k\GP$-mod. Then
$$
V_{\gp}(\m)=\biguplus_{(E,\C)}V^{\gp}_{E\times\C}(\m)^+,
$$
where the index runs over the set of isomorphism classes of Quillen
pairs. Moreover we have a homeomorphism
$$
V_{E,\C}(\m)^+/W_G(E,\C)\simeq V^{\gp}_{E\times\C}(\m)^+.
$$

\begin{proof} We know for $\m=\k$ the Quillen stratification
$$
V_{\gp}=\biguplus_{(E,\C)}{V^{\gp}_{E\times\C}}^+,
$$
with
$$
V_{E\times\C}^+/W_G(E,\C)\simeq {V^{\gp}_{E\times\C}}^+.
$$
From Remark 5.1.3 we have
$$
V_{\gp}(\m)=\bigcup_{(E,\C)}\iota^{\gp}_{E\times\C}V_{E\times\C}(\m).
$$
Since Lemma 5.1.4 implies
$$
V_{E\times\C}(\m)^+=V_{E\times\C}(\m)-\bigcup_{E'\subsetneqq
E}\iota^E_{E'}V_{E'\times(\C|E')}(\m),
$$
we get
$$
V_{\gp}(\m)=\bigcup_{(E,\C)}\iota^{\gp}_{E\times\C}V_{E\times\C}(\m)^+.
$$
It follows from Quillen's original result that, if $(E,\C)$ runs
over the set of isomorphism classes of Quillen pairs,
$$
V_{\gp}(\m)=\biguplus_{(E,\C)}V^{\gp}_{E\times\C}(\m)^+,
$$
and that there is a homeomorphism
$$
V_{E,\C}(\m)^+/W_G(E,\C)\simeq V^{\gp}_{E\times\C}(\m)^+.
$$
\end{proof}
\end{theorem}

\subsection{Subcategories and tensor products}

Here we record a few more expected interesting results as
consequences of the Quillen stratification.

\begin{theorem} Suppose $\P$ is a $G$-poset and $H$ is a subgroup
of $G$. Then, for any $\m\in k\GP$-mod,
$$
{\iota^{\gp}_{H\propto\P}}^{-1}V_{\gp}(\m)=V_{H\propto\P}(\m).
$$

\begin{proof} By Corollary 4.3.3 $\iota^{\gp}_{H\propto\P}V_{H\propto\P}(\m)\subset V_{\gp}(\m)$.
Given a Quillen pair $(E,\C)$ for $H\propto\P$, $(E,\C)$ is also a
Quillen pair for $\gp$ with $W_H(E,\C)\subset W_G(E,\C)$. By Theorem
5.1.5, there is a commutative diagram
$$
\xymatrix{V_{E\times\C}(\m)^+/W_H(E,\C) \ar@{>>}[rr] \ar[d]_{\simeq} &&
V_{E\times\C}(\m)^+/W_G(E,\C) \ar[d]^{\simeq}\\
V^{H\propto\P}_{E\times\C}(\m)^+ \ar[rr]_{\iota^{\gp}_{H\propto\P}} &&
V^{\gp}_{E\times\C}(\m)^+}
$$
Since the upper horizontal map is surjective and the vertical maps are homeomorphisms, the lower
horizontal map is surjective. We are done.
\end{proof}
\end{theorem}

We could not find a way to generalize the above statement for an
arbitrary sub-transporter category of $\gp$. However we have the
following.

\begin{proposition} Suppose $\P$ is a $G$-poset and $\Q$ is a $G$-subposet. Then, for any
$M\in kG$-mod,
$$
{\iota^{\gp}_{G\propto\Q}}^{-1}V_{\gp}(\kappa_M)=V_{G\propto\Q}(\kappa_M).
$$

\begin{proof} The proof is similar to the previous one. If $(E,\C)$ is a Quillen
pair for $G\propto\Q$, then there exists a unique Quillen pair
$(E,\D)$ for $\gp$. We note that $V_{E\times\C}(\kappa_M)^+
=V_{E\times\D}(\kappa_M)^+$. Also by definition if $g\C=\C$ then
$g\D=\D$ as well, for any $g\in G$. It implies that
$W_G(E,\C)\subset W_G(E,\D)$.
\end{proof}
\end{proposition}

Recall that for any category $\D$ there is a diagonal functor
$\triangle : \D\to\D\times\D$. Then image of $\D$ is written as
$\triangle\D$. If $\C$ is a connected poset, then so is
$\C\times\C$. Suppose $\gp$ is a transporter category and $(E,\C)$
is a Quillen pair for the $G$-poset $\P$. Then $(\triangle
E,\C\times\C)$ is a Quillen pair for the $(G\times G)$-poset
$\P\times\P$. Note that from Remark 2.1.2 $(G\times
G)\propto(\P\times\P)\cong\GP\times\GP$.

\begin{corollary} Suppose $G$ is a finite group and $\P$ is a finite $G$-poset.
If $H$ is a subgroup of $G$ and $\Q$ is a $H$-subposet of $\P$, then
${\iota^{\gp}_{H\propto\Q}}^{-1}V_{\gp}(\kappa_M)=V_{H\propto\Q}(\kappa_M)$.
Particularly for any $k(G\times G)$-module $N$
$$
{\iota^{(G\times G)\propto(\P\times\P)}_{\gp}}^{-1}V_{(G\times
G)\propto(\P\times\P)}(\kappa_N)=V_{\gp}(\kappa_N).
$$

\begin{proof} We apply Theorem 5.2.1 and Proposition 5.2.2 to the following three transporter
categories $H\propto\Q\subset H\propto\P\subset\gp$.

For the special case, we consider $\gp\cong\triangle
G\propto\triangle\P\subset\triangle
G\propto(\P\times\P)\subset(G\times G)\propto(\P\times\P)$.
\end{proof}
\end{corollary}

We need a technical result for the last main result.

\begin{proposition} Suppose $\D$ and $\C$ are two finite categories
with finitely generated ordinary cohomology rings. If $\m, \m' \in
k\D$-mod and $\n, \n' \in k\C$-mod, then
$$
\Ext^*_{k(\D\times\C)}(\m\otimes\n,\m'\otimes_k\n')\cong\Ext^*_{k\D}(\m,\m')
\otimes_k\Ext^*_{k\C}(\n,\n').
$$

\begin{proof} Let $\p^{\D}_* \to \m \to 0$ and $\p^{\C}_*\to\n\to 0$
be two projective resolutions. Then $\p^{\D}_*\otimes_k\p^{\C}_* \to
\m\otimes_k\n \to 0$ is a projective resolution of the
$k(\D\times\C)=k\D\otimes_kk\C$-module $\m\otimes_k\n$. The
isomorphism follows from
$$
\Hom_{k\D\otimes k\C}(\p^{\D}_*\otimes_k\p^{\C}_*,\m'\otimes_k
\n')\cong\Hom_{k\D}(\p^{\D}_*,\m')\otimes_k\Hom_{k\C}(\p^{\C}_*,\n'),
$$
and the K\"unneth formula.
\end{proof}
\end{proposition}

When we examine transporter categories, the above result has the
following consequences.

\begin{corollary} Let $G_1\propto\P_1$ and $G_2\propto\P_2$ be two
transporter categories. Then there is a natural isomorphism
$$
V_{(G_1\times
G_2)\propto(\P_1\times\P_2)}=V_{(G_1\propto\P_1)\times(G_2\propto\P_2)}
=V_{G_1\propto\P_1}\times V_{G_2\propto\P_2}.
$$
If $M\in kG_1$-mod and $N\in kG_2$-mod, under the above isomorphism
we have furthermore
$$
V_{(G_1\times G_2)\propto(\P_1\times\P_2)}(\kappa_{M\otimes
N})=V_{(G_1\propto\P_1)\times(G_2\propto\P_2)}(\kappa_M\otimes_k\kappa_N)
=V_{G_1\propto\P_1}(\kappa_M)\times V_{G_2\propto\P_2}(\kappa_N).
$$

\begin{proof} From Remark 2.1.2 we have $(G_1\times
G_2)\propto(\P_1\times\P_2)\cong(G_1\propto\P_1)\times(G_2\propto\P_2)$.
The statements on varieties are true, because
$$
\begin{array}{ll}
\Ext^*_{k[(G_1\times G_2)\propto(\P_1\times\P_2)]}(\kappa_{M\otimes
N},\kappa_{M\otimes N})
&\cong\Ext^*_{k[(G_1\propto\P_1)\times(G_2\propto\P_2)]}(\kappa_M\otimes\kappa_N,\kappa_M\otimes\kappa_N)\\
&\cong\Ext^*_{k(G_1\propto\P_1)}(\kappa_M,\kappa_M)\otimes\Ext^*_{k(G_2\propto\P_2)}(\kappa_N,\kappa_N).
\end{array}
$$
\end{proof}
\end{corollary}

Based on the preceding results, we can extend the tensor product
formula $V_G(M\otimes_k N)=V_G(M)\cap V_G(N)$ of Avrunin-Scott
\cite{AS}.

\begin{theorem} Let $G$ be a finite group and $\P$ a finite $G$-poset. Suppose
$M, N$ are two $kG$-modules. Then
$V_{\gp}(\kappa_M\hotimes\kappa_N)=V_{\gp}(\kappa_M)\cap
V_{\gp}(\kappa_N)$.

\begin{proof} Let us consider the functor $\triangle : \gp \to \GP\times\GP$
and the restriction along it. Because
$$
\kappa_M\hotimes\kappa_N =
{\Res}^{\GP\times\GP}_{\gp}(\kappa_M\otimes\kappa_N),
$$
by Corollary 5.2.3 we have
$$
V_{\gp}(\kappa_M\hotimes_k\kappa_N)={\iota^{\GP\times\GP}_{\gp}}^{-1}V_{\GP\times
\GP}(\kappa_M\otimes_k\kappa_N).
$$
From $(G\times G)\propto(\P\times\P)\cong\GP\times\GP$ and Corollary
5.2.5 the right-hand-side is
$$
{\iota^{\GP\times\GP}_{\gp}}^{-1}\{V_{\gp}(\kappa_M)\times
V_{\gp}(\kappa_N)\},
$$
which exactly is $V_{\gp}(\kappa_M)\cap V_{\gp}(\kappa_N)$.
\end{proof}
\end{theorem}

With Theorem 5.2.6, taking Theorem 4.2.1 (5) in account one can see
there is some relation between the conditions in (9a) and (9b) in
Theorem 4.2.1.

\end{document}